\documentclass[hidelinks,onefignum,onetabnum]{siamart220329}

\def\dfo#1#2{\frac{\mathrm{d} {#1}}{\mathrm{d} {#2}}}
\def\dft#1#2{\frac{\mathrm{d} ^2{#1}}{\mathrm{d} {#2} ^2}}
\def\dfn#1#2#3{\frac{\mathrm{d}^{#3}{#1}}{\mathrm{d}{#2}^{#3}}}
\def\pfo#1#2{\frac{\partial {#1}}{\partial {#2}}}
\def\pft#1#2{\frac{\partial^2{#1}}{\partial{#2} ^2}}

\def\mbf#1{\mathbf{#1}}
\def\est{\mathrm{e.s.t.}}
\def\updown#1#2#3{{#1}^{(#3)}_{#2}}  
\def\fc{g_*'}
\def\fcc{g_*''}
\def\fccc{g_*'''}

\allowdisplaybreaks

\usepackage{lipsum}
\usepackage{amsfonts}
\usepackage{graphicx}
\usepackage{epstopdf}
\usepackage{algorithmic}
\usepackage{IEEEtrantools}
\ifpdf
  \DeclareGraphicsExtensions{.eps,.pdf,.png,.jpg}
\else
  \DeclareGraphicsExtensions{.eps}
\fi

\newsiamremark{remark}{Remark}
\newsiamremark{hypothesis}{Hypothesis}
\crefname{hypothesis}{Hypothesis}{Hypotheses}
\newsiamthm{claim}{Claim}

\headers{Weakly Nonlinear Theory for Patterning in Structured Models}{Ridgway, Dalwadi, Pearce, Chapman}

\title{A Weakly Nonlinear Theory for Pattern Formation in Structured Models with Localized Solutions}

\author{Wesley J. M. Ridgway\thanks{Department of Mathematics and Statistics, University of Saskatchewan}
\and Mohit P. Dalwadi\thanks{Mathematical Institute, University of Oxford, Oxford, OX2 6GG, UK}
\and {Philip Pearce\thanks{D\MakeLowercase{epartment of} M\MakeLowercase{athematics}, U\MakeLowercase{niversity} C\MakeLowercase{ollege} L\MakeLowercase{ondon}, WC1H 0AY, UK}} \thanks{Institute for the Physics of Living Systems, University College London, London, UK} \and S. Jonathan Chapman\footnotemark[2]}

\ifpdf
\hypersetup{
  pdftitle={weakly nonlinear},
  pdfauthor={Wesley J. M. Ridgway, Mohit P. Dalwadi, Philip Pearce, S. Jonathan Chapman}
}
\fi

\begin{document}

\maketitle

\begin{abstract}

   Structured models, such as PDEs structured by age or phenotype, provide a setting to study pattern formation in heterogeneous populations. Classical tools to quantify the emergence of patterns, such as linear and weakly nonlinear analyses, pose significant mathematical challenges for these models due to sharply peaked or singular steady states. Here, we present a weakly nonlinear framework that extends classical tools to structured PDE models in settings where the base state is spatially uniform, but exponentially localized in the structured variable. Our approach utilizes WKBJ asymptotics and an analysis of the Stokes phenomenon to systematically resolve the solution structure in the limit where the steady state tends to a Dirac-delta function. To demonstrate our method, we consider a chemically structured (nonlocal) model of motile bacteria that interact through quorum sensing. For this example, our analysis yields an amplitude equation that governs the solution dynamics near a linear instability, and predicts a pitchfork bifurcation. From the amplitude equation, we deduce an effective parameter grouping whose sign determines whether the pitchfork bifurcation is subcritical or supercritical. Although we demonstrate our framework for a specific example, our techniques are broadly applicable.

\end{abstract}


\section{Introduction}

Heterogeneities within a population, such as differences in age, size and shape, affect individual interactions and population-level behaviors. For example, phenotypic heterogeneity is thought to be an important aspect of therapeutic resistance in cancer cells \cite{marusyk2012intra}, and antibiotic resistance in bacteria \cite{dewachter2019bacterial,smith2023bacterial}.  Structured PDE models provide a framework to capture heterogeneity within a population by representing population densities in terms of an internal state, in addition to space and time \cite{kotbook,murraybook,auger2008structured}. Examples of structured models include age-structured models in demographics,  epidemiology, and cell biology \cite{m1925applications,murraybook,inaba2017age,perthame2007transport}, size-structured population models \cite{cushing1998introduction,perthame2007transport}, trait or phenotype-structured models of cell migration and selection-mutation systems \cite{lorenzi2020asymptotic,lorenzi2025phenotype,phan2024direct,diekmann2005dynamics,turanova2015model,lorenzi2022invasion,arnold2012existence,crossley2025modelling}, stemness-structured models of cancer \cite{celora2023spatio}, kinetic transport equations and chemotaxis models with `orientation' and phenotype-structuring \cite{loy2024hamilton,lorenzi2024phenotype}, and chemically structured models in quorum sensing systems \cite{ridgway2023motility,Hu2021}. Broadly, this previous work has shown that heterogeneity in the characteristics or states of individuals can have a significant impact on spatio-temporal pattern formation compared with homogeneous populations. However, analytic techniques for understanding the emergence of patterns in such problems are lacking, due to the technical challenges associated with the analysis of structured PDEs.

The onset of patterning in these systems often signifies physically or biologically relevant phenomena, such as phase transitions or population extinction.  Quantifying the onset of patterning is a classical problem that dates back to Alan Turing \cite{turing}, though significant progress has been made since then (see \cite{krause2021modern} for a recent survey). Common analyses of Turing systems involve detecting the emergence of spatio-temporal patterning from a spatially uniform base state, which gives a broad idea of a system's propensity to form patterns. Mathematically, characterisation of patterning can be accomplished through linear stability analysis and weakly nonlinear analysis; the former is used to detect bifurcations, and the latter reduces the local dynamics to a canonical amplitude equation \cite{strogatz1994nonlinear}. The primary obstacle in performing such analyses for structured models is that spatially uniform states can be non-uniform with respect to the internal state.  Moreover, many structured models admit solutions that tend to become sharply peaked or singular with respect to the internal state \cite{lorenzi2020asymptotic,barles2009concentration,diekmann2005dynamics,lorenzi2022trade,lorenzi2025derivation,lorz2011dirac}. Such localized solutions represent homogeneous states in the sense that all individuals have nearly identical characteristics. Standard methods for linear and weakly nonlinear analyses often cannot be applied in such settings because the linearised system will have non-constant coefficients. While analytic techniques to deal with spatial and spatio-temporally evolving coefficients have been developed for local PDEs \cite{krause2024pattern, krause2020one,dalwadi2023universal}, many structured models are also nonlocal, adding another layer of difficulty.  As such, new mathematical tools are needed to extend classical Turing stability analyses to base states that are spatially uniform but non-constant in the structured variable. 

Previous work in this area includes linear stability analyses in situations where the base state is constant in both space and internal state \cite{lorenzi2025pattern,genieys2006adaptive,perthame2007concentration}. In non-structured models, sharply peaked  ``spike'' solutions can be analysed with a well-developed weakly nonlinear theory \cite{gomez2021hopf, wong2020weakly, kolokolnikov2021competition, veerman2015breathing}; example applications include various reaction-diffusion systems \cite{iron2001stability,wei2002spikes,wei2001spikes,chen2011stability,doelman1997pattern,doelman1998stability, al2022spikes,kolokolnikov2009spot, jabin2023collective} and kinetic transport models \cite{loy2024hamilton}. In structured models, previous studies have also proved global stability of Dirac-delta function solutions of well-mixed systems \cite{lorenzi2020asymptotic,barles2009concentration} as well as proved the existence of patterned states \cite{arnold2012existence,turanova2015model} and bifurcations from constant base states \cite{cushing1985equilibria}. The linear stability and bifurcation structure of non-constant steady states for a class of age-structured models with spatial diffusion was tackled with an abstract semigroup approach  (see \cite{walker2025stability,walker2010global} and references therein), where the stability is quantified in terms of the spectrum of a derived nonlocal eigenvalue problem.  Additionally, hybrid analytic-numeric methods have been developed to study traveling-wave patterns \cite{lorenzi2022trade,lorenzi2025derivation,benichou2012front}. Despite this progress, there is no general framework for explicitly quantifying pattern formation in structured models when the base state is sharply peaked or singular in the internal state.

Here, we present a perturbative framework for linear and weakly nonlinear analysis of structured population models in settings where the steady state is exponentially localized in the structured variable, but spatially uniform.  We demonstrate the applicability of our framework for a chemically structured model of motile bacteria that interact through quorum sensing. In the limit of a certain small parameter, the steady states in this model tend to a Dirac delta-function with respect to the internal state, but are uniform in space.  Our analysis overcomes the challenge of a non-constant, singular base state through the use of WKBJ asymptotics and an analysis of the Stokes phenomenon, which systematically resolve the localized structure of the internal state. This analysis yields a criterion for the onset of patterning  and quantifies the local dynamics near the bifurcation.  Although we consider a specific example, our framework is broadly applicable, particularly in settings where previous methods would not be appropriate.

The structure of the rest of this paper is as follows. In Section \ref{sec:model}, we present the model of motile quorum sensing bacteria, which we use to demonstrate our theory. In Section \ref{sec:prelim} we construct the exponentially localized steady state and   show that it formally tends to a Dirac delta in the limit of a certain small parameter. We then determine the location of the bifurcation through a linear stability analysis. Our main analysis is presented in Section \ref{sec:weakly_nonlinear}, where we perform a weakly nonlinear analysis near the bifurcation. The result of this analysis is an amplitude equation governing the spatio-temporal dynamics near the instability. The amplitude equation is that of a canonical pitchfork bifurcation, whose coefficients dictate whether the bifurcation is supercritical or subcritical. In Section \ref{sec:numerical_results}, we verify our predictions with numerical solutions of the governing equations and interpret our results in a biophysical context. Finally, in Section \ref{sec:discussion}, we discuss generalisations of our method as well as its range of applicability.

\section{An example model} \label{sec:model}

To demonstrate our weakly nonlinear theory, it is instructive to consider a specific example, though we emphasize that our framework is applicable to a wider class of models. For our example, we use the chemically structured population model of signalling bacteria in \cite{ridgway2023motility, ridgwaythesis}. This model  physically describes a population of cells whose motility is effectively regulated by quorum sensing (QS). The cells generate signalling molecules, commonly called autoinducers (AIs). This chemical signal affects the intracellular gene-regulatory kinetics, which in turn regulates the motility.  We refer the reader to \cite{ridgway2023motility} for a more detailed discussion of the biology and model derivation \emph{via} formal upscaling from cell-level rules.

In one spatial dimension, the (dimensionless) governing equations read 
\begin{subequations}
\begin{align}
	\pfo{ n}{t} &= \mathcal{D}(u)\pft{n}{x} + \varepsilon\pft{n}{u} - \pfo{}{u}\big(f(u,c)n\big), \quad x\in\Omega, \quad u>0, \quad t>0, \label{eqn:gov_eq_n} \\
    \pfo{c}{t} &= D_c\pft{c}{x} - \beta c + \alpha_0\int_0^\infty u\, n(x,u,t)\,\mathrm{d}u, \quad x\in\Omega,\quad t >0,\label{eqn:dim_gov_eq_c}
\end{align}
\label{eqn:gov_eqs}%
\end{subequations}
where the independent variables $x$ and $t$ represent position and time. The cell density $n(x,u,t)$ is structured in terms of the intracellular concentration $u$ of an internal chemical (e.g.~a transcription factor or gene-regulatory protein). The diffusion coefficient $\mathcal{D}(u)>0$ describes the effect of gene-regulated ($u$-dependent) cellular motility. The function $f$ describes the intracellular gene-regulatory reaction kinetics. The parameter $0 < \varepsilon \ll 1$ is a state-space diffusion coefficient, a deterministic representation of a small stochastic component in the reaction kinetics. As we will later see, this term can be thought of as a regularization of \eqref{eqn:gov_eqs}. An external chemical signal with concentration $c(x,t)$ permeates the population. This chemical has a natural decay rate $\beta$, diffusion coefficient $D_c$, and is secreted by the population at a  rate $\alpha_0 u$. Since $u$ is the internal chemical concentration, the combined effect of the secretion over all $u$ generates the nonlocal integral term in \eqref{eqn:dim_gov_eq_c}. Finally, the population is assumed to be constrained to a finite domain $\Omega:=(0,L)$. We impose no-flux boundary conditions in both physical space $(x)$ and state space ($u$) as
\begin{subequations}
    \begin{align}
         \mathcal{D}(u)\pfo{n}{x} = D_c\pfo{c}{x} &= 0, \quad x = 0,\, L, \label{eqn:gov_BCs_x}\\
        \varepsilon\pfo{n}{u} - f(u,c)n &\to 0, \quad u\to 0, \,\infty. \label{eqn:gov_BCs_u}
    \end{align}
    \label{eqn:gov_BCs}%
\end{subequations}
where \eqref{eqn:gov_BCs_u} ensures physically realistic (bounded and positive) internal concentrations. For later convenience, we define the cell density $\rho(x,t)$ as 
\begin{equation}
    \rho(x,t):=\int_0^\infty n(x,u,t)\,\mathrm{d}u.
    \label{eqn:rho_def}
\end{equation}
This quantity effectively sums the full distribution of internal states $n(x,u,t)$ into a local cell density at each point in space, and will be useful in Section \ref{sec:numerical_results} where we construct bifurcation diagrams. Since the RHS of \eqref{eqn:gov_eq_n} is a total divergence in $(x,u)$-space, the global cell density is conserved. We denote the global cell density by
\begin{equation}
    \rho^*:= \frac{1}{L}\int_\Omega \rho(x,t)\,\mathrm{d}x.
    \label{eqn:normalisation}
\end{equation}

Finally, it is helpful to impose a functional form for the gene-regulatory kinetics $f$. To illustrate the theory, we assume the following form 
\begin{equation}
    f(u,c) = g(c) -\lambda u,
    \label{eqn:gov_kinetics}
\end{equation}
which is linear in $u$ but can be nonlinear in $c$ for generality. Here, $\lambda$ represents a natural decay rate and the function $g$ describes production of internal chemical. Many quorum sensing systems exhibit positive feedback \cite{dalwadi2021emergent}. Since an increase in $u$ corresponds to an increase in $c$ \emph{via} the integral term in \eqref{eqn:dim_gov_eq_c}, we can incorporate positive feedback by requiring that an increase in $c$  corresponds to an increase in $u$. This can be encoded by requiring the production function $g$ to be increasing, i.e.~$g'(c)>0$ for all $c\geq0$. To ensure non-negative and finite  concentrations, we  assume further that $g$ is positive and bounded. Our full model consists of Eqs.~\eqref{eqn:gov_eqs}--\eqref{eqn:gov_kinetics}.

Our remaining goal is to demonstrate our weakly non-linear theory by analysing the solution structure near a linear instability in the spatially uniform steady state of Eqs.~\eqref{eqn:gov_eqs}--\eqref{eqn:gov_BCs}. In \cite{ridgway2023motility} it was shown that the spatially uniform solution can become unstable when $\mathcal{D}'$ is sufficiently large and negative, and that the dynamics following the instability lead to motility-induced phase separation (MIPS) \cite{cates2015motility}. A linear stability theory was developed to quantify the instability condition and the bifurcating branches of steady states computed numerically.  Here, we extend the analysis to include non-linear terms, thereby allowing us to characterise the bifurcation  as a symmetric pitchfork, as well as derive an explicit formula that dictates whether the bifurcation is supercritical or subcritical. We also explicitly calculate the local behaviour of the bifurcating solution branches. Our analysis also generalises the linear theory in \cite{ridgway2023motility} to GRN kinetics with more a generic positive feedback term $g(c)$ in \eqref{eqn:gov_kinetics}. 

Before getting into the details of our analysis, we give a brief outline of our methodology.  The governing equations \eqref{eqn:gov_eqs}--\eqref{eqn:gov_BCs} admit a $u$-dependent, spatially uniform solution which formally tends to a Dirac-delta in the limit $\varepsilon\to0^+$. For small finite $\varepsilon$, this solution is instead exponentially localized in $u$ to an $\mathcal{O}(\varepsilon^{1/2})$ `inner' region around a point $u^*$ to be determine as part of the analysis. To study the nonlinear stability of this solution, we introduce small perturbations near the instability. For the perturbation to remain small compared to the base state, each term in the perturbation expansion must also be exponentially localized. We therefore seek a WKBJ solution for each term in the perturbation expansion, whereby the WKBJ amplitudes satisfy a set of singularly perturbed ODEs. Through a careful analysis of the Stokes switching in the outer solution for the amplitudes, we exchange the boundary conditions in the outer region for a regularity condition in the inner region. This allows us to construct power series solutions to the WKBJ amplitude equations in the inner region that satisfy the regularity condition. Then using the fact that the perturbations are exponentially localized, we use Laplace's method to evaluate all integrals arising from the nonlocal term in \eqref{eqn:dim_gov_eq_c}. Our methodology  relies on the exponential localisation of the base state; it allows the nonlocal terms to be evaluated explicitly using only the local behaviour of the WKBJ amplitudes. 

\pagebreak

\section{Preliminary analysis}
\label{sec:prelim}

\subsection{Spatially homogeneous steady state}

The governing equations \eqref{eqn:gov_eqs}--\eqref{eqn:gov_BCs} admit a spatially homogeneous steady state $n(x,u,t)=n^*(u)$ given by 
\begin{subequations}
    \begin{align}
        n^*(u) &=N\exp\left[-\frac{\lambda}{2\varepsilon}(u - u^*)^2\right].
        \label{eqn:canonical_density_ss}
    \end{align}
The normalisation condition \eqref{eqn:normalisation} implies
\begin{equation}
    N = \rho^*\sqrt{\frac{\lambda}{2\pi\varepsilon}} + \est,
    \label{eqn:canonical_density_ss_normalisation}
\end{equation}
\label{eqn:n_ss}%
\end{subequations}
where here and henceforth $\est$~denotes exponentially small terms. We emphasize that the steady state \eqref{eqn:n_ss} formally tends to a Dirac delta-function in the limit $\varepsilon\to0^+$. Physically, this corresponds to the deterministic limit of the cell-level reaction kinetics. 

The steady AI concentration $c^*$ and mean internal concentration $u^*$ satisfy
\begin{subequations}
\begin{align}
    f(u^*,c^*)&:=g(c^*)- \lambda u^*=0, \label{eqn:steady_kinetics}\\
    \beta c^* &= \alpha_0u^*\rho^* + \est, \label{eqn:steady_c}
\end{align}
\label{eqn:steady_algebraic}
\end{subequations}
Depending on the functional form for $g$ in \eqref{eqn:gov_kinetics}, the algebraic equation 
\begin{equation}
    g(c^*) = \frac{\lambda\beta}{\alpha_0\rho^*}c^*,
    \label{eqn:c_ss}
\end{equation}
may have multiple solutions. We note that since $g$ is monotonically increasing by assumption, there will always be a parameter regime for $\rho^*$ where there is at least one positive solution.  We focus our attention on this parameter regime. The spatially uniform steady states of Eqs.~\eqref{eqn:gov_eqs}--\eqref{eqn:gov_BCs} are then given by positive solutions of \eqref{eqn:c_ss} for $c^*$ and \eqref{eqn:canonical_density_ss} with $u^*\sim g(c^*)/\lambda$ for $n^*(u)$, up to exponentially small corrections.

\subsection{Linear stability}

To perform the weakly nonlinear analysis near the bifurcation, we must first determine the location of the bifurcation. As demonstrated in \cite{ridgway2023motility}, a linear stability analysis around the spatially uniform steady states in \eqref{eqn:n_ss}, \eqref{eqn:c_ss} leads to the following algebraic equation for the eigenvalue $\sigma$: 
\begin{equation}
	\sigma + D_ck^2 + \beta - \alpha_0\fc\rho^*\frac{\sigma + \left(\mathcal{D}_{0*} - u^*\mathcal{D}_*'\right)k^2}{\left(\sigma + \mathcal{D}_{0*}k^2\right)\left(\sigma + \mathcal{D}_{0*}k^2 + \lambda\right)}=\mathcal{O}(\varepsilon),
	\label{eqn:eigenvalue_problem}
\end{equation}
where $\mathcal{D}_{*}:=\mathcal{D}(u^*)$, $\mathcal{D}_*':=\mathcal{D}'(u^*)$, $g_*':=g'(c^*)$, and wavenumber $k$ with allowed values equal to integer multiples of $\pi/L$ to satisfy the no-flux conditions \eqref{eqn:gov_BCs_x}. We hereafter use subscript `$*$' to denote evaluation at steady state, i.e.~at $u=u^*$ and $c=c^*$. Eq.~\eqref{eqn:eigenvalue_problem} is a generalisation of Eq.~(13) in \cite{ridgway2023motility} to to a generic production term $g(c)$. Since we are interested in instabilities that lead to spatial patterning, we focus on situations where the spatially uniform solution is stable to uniform perturbations. That is, we take $g(c)$ such that there are no unstable eigenvalues when $k=0$. This implies the following restriction on $g'(c^*)$
\begin{equation}
    \fc < \frac{\lambda \beta}{\alpha_0 \rho^*},
    \label{eqn:gdash_constraint}
\end{equation}
in addition to $\fc>0$ as previously assumed. 

For spatially non-uniform perturbations (i.e.~$k^2>0$), \eqref{eqn:eigenvalue_problem} admits a real-valued, positive eigenvalue $\sigma$ whenever $\mathcal{D}(u)$ is such that $\mathcal{D}_*$ and $\mathcal{D}'_*$ satisfy
\begin{equation}
	\mathcal{D}'_{*}< -\frac{\mathcal{D}_{*}}{u^*}\left[\frac{\left(D_ck^2 + \beta\right)\left(\mathcal{D}_{*}k^2 + \lambda\right)}{\alpha_0\rho^*g'_*} - 1\right] + \mathcal{O}(\varepsilon),
	\label{eqn:zero_eig}
\end{equation}
where now $k=\pi/L$ is the critical wavenumber for instability. Since the $\mathcal{O}(1)$ term on the RHS of \eqref{eqn:zero_eig} is negative (which follows from \eqref{eqn:gdash_constraint}), this bifurcation can only occur when $\mathcal{D}_*'$ is negative. Physically this means that spatio-temporal patterning emerges from a spatially uniform state when motility is sufficiently repressed by quorum sensing, as measured by the derivative of $\mathcal{D}$ at steady state.

To keep our analysis general, we work with a general $\mathcal{D}$ that depends smoothly on a bifurcation parameter $\theta$. We denote this parameter dependence by $\mathcal{D}=\mathcal{D}(u;\theta)$, where the bifurcation occurs at $\theta=0$ (without approximation in $\varepsilon$). At the bifurcation, we introduce the notation
\begin{equation}
    \mathcal{D}(u;0) = \mathcal{D}_0(u). 
\end{equation}
For small $\varepsilon$ the bifurcation condition is given by
\begin{subequations}
\begin{equation}
	\mathcal{D}'_{0*}= -\frac{\mathcal{D}_{0*}}{u^*}\left[\frac{\left(D_ck^2 + \beta\right)\left(\mathcal{D}_{0*}k^2 + \lambda\right)}{\alpha_0\rho^*g'_*} - 1\right] + \mathcal{O}(\varepsilon),
    \label{eqn:bif_condition1}
\end{equation}
where $\mathcal{D}_{0*}$ and $\mathcal{D}_{0*}'$ are defined by  
\begin{equation}
    \mathcal{D}_0(u^*)=\mathcal{D}_{0*}, \quad \mathcal{D}_{0}'(u^*)=\mathcal{D}_{0*}'.
    \label{eqn:bif_condition2}
\end{equation}
\label{eqn:bif_condition}
\end{subequations}

\section{Weakly non-linear analysis near the bifurcation} \label{sec:weakly_nonlinear}

To study the local solution structure near the instability, we consider a perturbation to $\theta$ near the bifurcation. We therefore perturb $\theta$ as
\begin{equation}
    \theta = \theta_0\delta^2, 
    \label{eqn:theta_perturb}
\end{equation}
where $0<\delta\ll 1$ is a small parameter that will measure the size of the perturbation to the steady state. We absorb the parameter $\theta_0$ into a function $d(u)$ by Taylor expanding $\mathcal{D}$ near the bifurcation as
\begin{equation}
    \mathcal{D}(u;\theta_0\delta^2) = \mathcal{D}_{0}(u) + d(u)\delta^2 + \mathcal{O}(\delta^4), \quad d(u):=\theta_0\pfo{\mathcal{D}}{\theta}(u;0). 
    \label{eqn:WNA_Dp_perturb}
\end{equation}
For later convenience, we will assume that if $\theta_0\neq0$, then
\begin{equation}
    d(u^*)=0, \quad d'(u^*)\neq0,
    \label{eqn:d_perturb}
\end{equation}
so that at $u=u^*$, a small finite $\delta$ perturbs the derivatives of $\mathcal{D}$ at $\mathcal{O}(\delta^2)$, but introduces only a higher order $\mathcal{O}(\delta^4)$ perturbation to the value of $\mathcal{D}$. The first assumption in \eqref{eqn:d_perturb} simplifies our calculations, but is not a technical necessity and does not fundamentally change the bifurcation structure. In Section \ref{sec:Odelta3}, we flag where this assumption comes into our analysis. The second assumption in \eqref{eqn:d_perturb} is a transversality condition ensuring that the bifurcation is crossed at non-zero speed as $\theta$ varies.

The choice of $\mathcal{O}(\delta^2)$ scaling in \eqref{eqn:theta_perturb} will be justified  in the course of our analysis. Briefly, we will find that there are no secular terms that appear at order $\mathcal{O}(\delta^2)$, which implies that any $\mathcal{O}(\delta)$ term in \eqref{eqn:theta_perturb} will yield a solvability condition that contains no information about the bifurcation structure beyond the linear theory. Since the system \eqref{eqn:gov_eqs} is invariant under the transformation   $x\mapsto L - x$, a symmetry argument implies that non-uniform solutions near the bifurcation come in pairs, so we expect a symmetric pitchfork bifurcation. 

In general, small perturbations around the uniform steady state will vary on an $\mathcal{O}(1)$ timescale. However, near a bifurcation the timescales involved are much longer due to the phenomenon of critical slowing \cite{strogatz1994nonlinear}. With the benefit of hindsight, the appropriate time rescaling is
\begin{equation}
	t = \delta^{-2} \tau,
	\label{eqn:WNA_t_scaling}
\end{equation}
where $\tau$ is slow time.

Next, we seek solutions of the governing equations \eqref{eqn:gov_eqs}--\eqref{eqn:gov_BCs} near the bifurcation. We therefore expand $n$ and $c$ in regular asymptotic series as
\begin{subequations}
\begin{align}
	n(x,u,\tau) &= n^*(u) + \delta\eta_1(x,u,\tau) + \delta^2\eta_2(x,u,\tau) + \delta^3\eta_3(x,u,\tau) + \mathcal{O}(\delta^4), \label{eqn:wna_expand_n}\\
	c(x,\tau) &= c^* + \delta c_1(x,\tau) + \delta^2 c_2(x,\tau) + \delta^3 c_3(x,\tau) + \mathcal{O}(\delta^4), \label{eqn:wna_expand_c}
\end{align}
\label{eqn:wna_expand}%
\end{subequations}
where the spatially uniform solutions $n^*$ and $c^*$ are given in \eqref{eqn:n_ss}--\eqref{eqn:steady_algebraic} and $\eta_j$, $c_j$, $j=1,2,3$ are unknown functions to be determined. 

The next step is to insert the scalings \eqref{eqn:WNA_Dp_perturb}, \eqref{eqn:WNA_t_scaling} and expansions \eqref{eqn:wna_expand} into the governing equations and collect powers of $\delta$. In our analysis, we perform the expansion in $\delta$ without approximation in $\varepsilon$, we then analyse the problems at each order in $\delta$ in the limit $\varepsilon\ll 1$. The remainder of this section is divided into three parts. In Section \ref{sec:leading_orders}, we cover the analysis at $\mathcal{O}(1)$ and $\mathcal{O}(\delta)$. In Section \ref{sec:Odelta2}, we consider the $\mathcal{O}(\delta^2)$ problem, where we show that the scalings \eqref{eqn:WNA_Dp_perturb} and \eqref{eqn:WNA_t_scaling} are correct because they generate a trivially satisfied solvability condition. Our weakly nonlinear analysis is concluded in Section \ref{sec:Odelta3} with an analysis of the terms at $\mathcal{O}(\delta^3)$, resulting in an amplitude equation for the $\tau$-dependence of $\eta_1$ and $c_1$ in \eqref{eqn:wna_expand}. This amplitude equation will be in the normal form for a symmetric pitchfork bifurcation. 

\subsection{Analysis at $\mathcal{O}(1)$ and $\mathcal{O}(\delta)$} \label{sec:leading_orders}
We insert the scalings \eqref{eqn:WNA_Dp_perturb}, \eqref{eqn:WNA_t_scaling} and expansions \eqref{eqn:wna_expand} into the governing equations \eqref{eqn:gov_eqs}--\eqref{eqn:gov_BCs} and collect powers of $\delta$. The terms at $\mathcal{O}(1)$ vanish by construction since $n^*$ and $c^*$ are steady state solutions. At $\mathcal{O}(\delta)$, we obtain
\begin{subequations}
\begin{equation}
	\mathcal{L}\begin{pmatrix}
	    \eta_1 \\
        c_1
	\end{pmatrix}:=
    \begin{pmatrix}
    \mathcal{D}_0(u) \pft{\eta_1}{x}  + \varepsilon\pft{\eta_1}{u}  - \pfo{}{u}\left[f(u,c^*)\eta_1\right]- c_1 \fc\dfo{n^*}{u} \\
	D_c\pft{c_1}{x}-\beta c_1 +\alpha_0\int_0^\infty u \eta_1(x,u,t)\,\mathrm{d}u
    \end{pmatrix} = \mbf{0},
    \label{eqn:gov_eta1_c1_eqn}
\end{equation}
and the no-flux boundary conditions
\begin{align}
	\varepsilon\pfo{\eta_1}{u} - f(u,c^*)\eta_1 &\to0, \quad u\to0, \infty, \label{eqn:gov_eta1_BC}\\
	\pfo{\eta_1}{x} =\pfo{c_1}{x}&= 0, \quad x=0,L.
\end{align}
\label{eqn:gov_wna_Odelta}%
\end{subequations}
We seek non-trivial separable solutions to \eqref{eqn:gov_wna_Odelta} in the form
\begin{equation}
	\eta_1(x,u,\tau) = A(\tau) \eta(u)\cos(kx), \qquad c_1(x) = A(\tau) \cos(kx),
	\label{eqn:wna_Odelta_separation}
\end{equation}
where $\eta$ and $A$ are unknown functions to be determined. The function $A$, and its dependence on the local bifurcation parameter $d'_*$, determine the local solution structure near the bifurcation. Since $\tau$ is only a parameter in \eqref{eqn:gov_wna_Odelta}, $A(\tau)$ will be determined from our analysis of higher-order terms in $\delta$. Using the ansatz \eqref{eqn:wna_Odelta_separation} in \eqref{eqn:gov_wna_Odelta}, we obtain the following system for $\eta$:
\begin{subequations}
\begin{align}
	\varepsilon\dft{\eta}{u} + \dfo{}{u}\left[\lambda (u - u^*) \eta\right] - \mathcal{D}_0(u)k^2\eta &= \fc\dfo{n^*}{u}, \label{eqn:wna_eta_eqn} \\
	D_ck^2 + \beta  - \alpha_0 \int_0^\infty u\eta(u)\,\mathrm{d}u &=0,  \label{eqn:wna_eta_alg_eqn} \\
	\varepsilon\dfo{\eta}{u} - f(u,c^*)\eta &\to0, \quad u\to0, \infty. \label{eqn:wna_eta_eqn_bc}
\end{align}
\label{eqn:wna_Odelta_eqn}%
\end{subequations}
Eq.~\eqref{eqn:wna_Odelta_eqn} appears to be overdetermined, but the apparent `extra' equation \eqref{eqn:wna_eta_alg_eqn} actually determines the position of the bifurcation. To see this, we will obtain $\eta$ directly from \eqref{eqn:wna_eta_eqn} in the limit $\varepsilon\ll 1$ and insert the result into \eqref{eqn:wna_eta_alg_eqn}. After evaluating the integral, the resulting algebraic equation will be the bifurcation condition \eqref{eqn:bif_condition}. 

The unperturbed ($\delta=0$) solution \eqref{eqn:wna_expand_n} is exponentially localized with respect to $\varepsilon$ near $u=u^*$. In order for the perturbation $\delta\eta_1$ to remain small in \eqref{eqn:wna_expand_n} compared to the steady state $n^*$ when $0 < \delta \ll 1$, we require $\eta$ to also be exponentially localized in $\varepsilon$. We therefore seek WKBJ-type solutions of the form
\begin{equation}
	\eta(u)= \varepsilon^{-\frac{3}{2}}q(u)\exp\left[-\frac{\lambda(u-u^*)^2}{2\varepsilon}\right], \label{eqn:wna_Odelta_wkbj}
\end{equation}
for an unknown amplitude $q(u)$ to be determined. The $\varepsilon^{-3/2}$ prefactor is introduced for later convenience so that $q=\mathcal{O}(1)$ as $\varepsilon\to 0$. The WKBJ method has been employed in previous structured models to resolve sharply peaked solutions, where the phase satisfies a Hamilton-Jacobi equation \cite{diekmann2005dynamics, perthame2007concentration,lorenzi2022trade}. The key difference here is that the primary difficulty lies in determining the amplitude $q$, rather than the phase.  Continuing, we substitute our WKBJ ansatz \eqref{eqn:wna_Odelta_wkbj} into the governing equation for $\eta$ in \eqref{eqn:wna_eta_eqn} and find that $q$ satisfies the singularly perturbed equation
\begin{equation}
	\varepsilon \dft{q}{u} - \lambda (u - u^*)\dfo{q}{u} - \mathcal{D}_0k^2 q = -\rho^*\fc\sqrt{\frac{\lambda^3}{2\pi}}(u - u^*).
	\label{eqn:canonical_q_ODE}
\end{equation}
For boundary conditions, we require that $q$ be polynomially bounded as $u\to 0, \infty$ so that $\eta$ is exponentially localized. This guarantees that the boundary conditions  \eqref{eqn:wna_eta_eqn_bc} on $\eta$ are satisfied up to exponentially small terms.

We seek solutions of \eqref{eqn:canonical_q_ODE} in the form of a regular asymptotic series of the form
\begin{equation}
	q(u) \sim \sum_{j=0}^\infty q_j(u)\varepsilon^j,
	\label{eqn:canonical_q_expand}
\end{equation}
where each $q_j$ must also be polynomially bounded as $u\to0,\infty$. By substituting the expansion \eqref{eqn:canonical_q_expand} into the amplitude equation \eqref{eqn:canonical_q_ODE}  and equating coefficients of powers of $\varepsilon$, we obtain a recursive sequence of ODEs for the $q_j$'s given by
\begin{subequations}
\begin{align}
	\lambda (u - u^*)q_0' + \mathcal{D}_0(u)k^2q_0  &= \rho^*\fc\sqrt{\frac{\lambda^3}{2\pi}}(u - u^*), \label{eqn:canonical_q0_ode}\\
	\lambda (u - u^*)q_j' + \mathcal{D}_0(u)k^2q_j &= q_{j-1}'', \qquad j=1,\ldots. \label{eqn:canonical_qj_ode}
\end{align}
\label{eqn:canonical_q0_q1_ode}%
\end{subequations}

Since $u=u^*$ is a singular point of \eqref{eqn:canonical_q0_q1_ode}, the general solution for $q_j$ is non-smooth at $u^*$ for all $j$. This poses a potential issue since $q$ itself is smooth. However, as we show in Appendix \ref{app:stokes}, the non-smooth solutions are formally ruled out. This is because non-smooth solutions would generate a factorial-over-power divergent series in \eqref{eqn:canonical_q_expand}, which would switch on an exponentially large contribution to $q$ \emph{via} the Stokes phenomenon \cite{chapman1998exponential}. This exponentially large term would violate the boundary condition that $q$ is polynomially bounded as $u\to 0,\infty$. Therefore, we effectively enforce the boundary conditions on $q$ by demanding that $q_j$ is smooth for all $j$. 

Given the smoothness of $q_j$, we now seek solutions of \eqref{eqn:canonical_q0_q1_ode} in the form of a regular power series centred at $u=u^*$. We therefore expand $q_j$ as
\begin{equation}
	q_j(u) = \sum_{\ell=0}^\infty \updown{q}{j}{\ell} (u - u^*)^\ell,
	\label{eqn:canonical_qj_ps}
\end{equation}
and additionally Taylor expand $\mathcal{D}_0$ as
\begin{equation}
	\mathcal{D}_0(u) = \sum_{\ell =0}^\infty\frac{\mathcal{D}_{0*}^{(\ell)}}{\ell!}(u - u^*)^\ell, \quad \mathcal{D}_{0*}^{(\ell)}:=\dfn{}{u}{\ell}\bigg|_{u=u^*}\mathcal{D}_0(u).
	\label{eqn:canonical_D_taylor}
\end{equation}
Using only the local behaviour of $q$ near $u=u^*$, Laplace's method can be used to generate an asymptotic approximation of the integral in \eqref{eqn:wna_eta_alg_eqn}  containing arbitrarily many algebraic orders of $\varepsilon$.  Therefore, the local solution \eqref{eqn:canonical_qj_ps} is sufficient for our leading order (in $\varepsilon$) analysis. Using standard power series methods, we obtain the following recursion relations for the coefficients $\updown{q}{j}{\ell}$
\begin{subequations}
\begin{align}
	 \updown{q}{0}{0}&=0, \\
	\updown{q}{0}{\ell} &= \dfrac{1}{\mathcal{D}_{0*}k^2 + \ell\lambda}\left[\delta_{\ell1}\sqrt{\dfrac{\lambda^3}{2\pi}}\rho^*\fc - k^2\sum_{m=1}^\ell\frac{\mathcal{D}_{0*}^{(m)}}{m!}\updown{q}{0}{\ell - m}\right], \quad \ell = 1,\ldots, \label{eqn:wna_q0_ps_coeffs}\\
	\updown{q}{j}{0} &= \dfrac{2\updown{q}{j-1}{\ell+2}}{\mathcal{D}_{0*}k^2}, \quad j=1,\ldots,\\
	\updown{q}{j}{\ell}&=\dfrac{(\ell +1)(\ell + 2)\updown{q}{j-1}{\ell + 2} - k^2\sum\limits_{m=1}^\ell \updown{q}{j}{\ell - m}\dfrac{\mathcal{D}_*^{(m)}}{m!}}{\mathcal{D}_{0*}k^2 + \ell\lambda}, \quad j, \ell=1,\ldots, 
\end{align}
\label{eqn:canonical_wna_qj_coeffs}%
\end{subequations}
where $\delta_{ij}$ denotes the Kronecker delta.  The recursion relations in \eqref{eqn:canonical_wna_qj_coeffs} determine $\eta$ from \eqref{eqn:wna_Odelta_wkbj} locally near $u=u^*$. Inserting this result for $\eta$ into \eqref{eqn:wna_eta_alg_eqn} and utilising Laplace's method to evaluate the leading-order contribution to the integral, we recover \eqref{eqn:bif_condition}.  This concludes our analysis at linear order in $\delta$. As expected, this analysis does not determine $A(\tau)$ in \eqref{eqn:wna_Odelta_separation}. We therefore proceed to next order in $\delta$. 

\subsection{Analysis at $\mathcal{O}(\delta^2)$} \label{sec:Odelta2}

Next, we obtain governing equations for  $\eta_2$ and $c_2$ by inserting the expansion \eqref{eqn:wna_expand} and scalings \eqref{eqn:WNA_Dp_perturb}, \eqref{eqn:WNA_t_scaling} into the governing equations \eqref{eqn:gov_eqs}--\eqref{eqn:gov_BCs} and collecting $\mathcal{O}(\delta^2)$ terms. This yields
\begin{subequations}
\begin{equation}
	\mathcal{L}\begin{pmatrix}
	    \eta_2\\
        c_2
	\end{pmatrix} = 
    \begin{pmatrix}
        c_1^2\frac{\fcc}{2}\dfo{n^*}{u}  + c_1 \fc\pfo{\eta_1}{u}  \\
        0
    \end{pmatrix},
    \label{eqn:gov_canonical_wna_eta2_eqn}
\end{equation}
and the boundary conditions 
\begin{align}
	\varepsilon \pfo{\eta_2}{u} - f(u,c^*)\eta_2 &\to 0, \quad u\to0,\infty, \label{eqn:gov_canonical_eta2_BC}\\
	\pfo{\eta_2}{x} =\pfo{c_2}{x}&= 0, \quad x=0,L,
	\label{eqn:gov_canonical_c2_BC}
\end{align}
\label{eqn:gov_canonical_wna_Odelta2}%
\end{subequations}
where we have ignored exponentially small terms in \eqref{eqn:gov_canonical_eta2_BC}--\eqref{eqn:gov_canonical_c2_BC}. The homogeneous version of \eqref{eqn:gov_canonical_wna_Odelta2} has a non-trivial solution given by \eqref{eqn:wna_Odelta_separation}, where  $\eta$ in \eqref{eqn:wna_Odelta_separation} is given in \eqref{eqn:wna_Odelta_wkbj}, \eqref{eqn:canonical_q_expand}, \eqref{eqn:canonical_qj_ps}, and \eqref{eqn:canonical_wna_qj_coeffs}. The Fredholm Alternative therefore implies that \eqref{eqn:gov_canonical_wna_Odelta2} has no solution unless the RHS is orthogonal to the nullspace of the adjoint operator $\mathcal{L}^*$, which we now define. 

To obtain the adjoint operator, we first define an appropriate inner product. For two vector functions $\mbf{v}=(v_1(x,u), v_2(x))$, $\mbf{w}=(w_1(x,u),w_2(x))$ we define the inner product $\langle\mbf{v},\mbf{w}\rangle$ as
\begin{equation}
	\left\langle\mbf{v},\mbf{w}\right\rangle:= \int_{\Omega\times\mathbb{R}_+} v_1(x,u)w_1(x,u)\,\mathrm{d}u\mathrm{d}x + \int_{\Omega} v_2(x)w_2(x)\,\mathrm{d}x.
\end{equation}
Suppose that $\mbf{v}$ satisfies the no-flux boundary conditions
\begin{subequations}
\begin{align}
	\varepsilon \pfo{v_1}{u} - f(u,c^*)v_1 &\to 0, \quad u\to0,\infty, \label{eqn:v1_BC}\\
	\pfo{v_1}{x} =\pfo{v_2}{x}&= 0, \quad x=0,L.
	\label{eqn:v2_BC}
\end{align}
\label{eqn:v_BCs}%
\end{subequations}
Then the adjoint operator $\mathcal{L}^*$ is uniquely defined by $\langle\mathcal{L}(\mbf{v}),\mbf{w}\rangle = \langle\mbf{v},\mathcal{L}^*(\mbf{w})\rangle$ in combination with the adjoint boundary conditions 
\begin{subequations}
\begin{align}
	\exp\left(-\frac{\lambda(u - u^*)^2}{2\varepsilon}\right) \pfo{w_1}{u} &\to0, \quad u\to0, \infty, \label{eqn:WNA_V_adj_BC_u}\\
	\pfo{w_1}{x} =\pfo{w_2}{x}&= 0, \quad x=0,L. \label{eqn:WNA_V_adj_BC_x}
\end{align}
\label{eqn:WNA_V_adj_BC}%
\end{subequations}
It is sufficient to require that $w_1$ is polynomially bounded to satisfy  \eqref{eqn:WNA_V_adj_BC_u}. Integrating by parts in the definition of the adjoint operator, we find that $\mathcal{L}^*$ is given by
\begin{equation}
\mathcal{L}^*(\mbf{w}):=
\begin{pmatrix}
	\mathcal{D}_0\pft{w_1}{x}  +\varepsilon \pft{w_1}{u} + f(u,c^*)\pfo{w_1}{u}  + \alpha_0 u w_2 \\
	D_c\pft{w_2}{x} - \beta w_2 - \int_0^\infty \fc\dfo{n^*}{u}w_1(x,u)\,\mathrm{d}u
    \end{pmatrix}.
	\label{eqn:WNA_adj_op}
\end{equation}

In order to explicitly apply the Fredholm Alternative to \eqref{eqn:gov_canonical_wna_eta2_eqn}, we must obtain all non-trivial solutions of $\mathcal{L}^*(\mbf{w})=\mbf{0}$. One such solution is given by $(w_1, w_2) =(1,0)$ up to an exponentially small correction. This `trivial' solution is present because the governing equation \eqref{eqn:gov_eq_n} is a conservation equation; it has no relevance to the bifurcation and yields a trivial solvability condition as we will later see. To obtain any remaining solutions of $\mathcal{L}^*(\mbf{w})=\mbf{0}$, we seek Fourier mode solutions in the form 
\begin{equation}
	w_1(x,u) = W(u)\cos(k x), \qquad w_2(x) = \cos(k x),
	\label{eqn:wna_adj_mode}
\end{equation}
where we have fixed the scaling of $w_2$ for convenience. Next, we insert \eqref{eqn:wna_adj_mode} into $\mathcal{L}^*(\mbf{w})=\mbf{0}$ and \eqref{eqn:WNA_adj_op} and find that $W$ satisfies
\begin{subequations}
\begin{align}
	\varepsilon \dft{W}{u} - \lambda (u - u^*) \dfo{W}{u}  - \mathcal{D}_0(u)k^2 W &= -\alpha_0 u, \label{eqn:wna_adj_W_eqn}\\
	D_c k^2 + \beta  + \int_0^\infty \fc \dfo{n^*}{u} W(u) \,\mathrm{d}u &=0,	\label{eqn:wna_adj_bif_condition}
\end{align}
\label{eqn:wna_adj_W_problem}%
\end{subequations}
with $W$ polynomially bounded as $u\to 0,\infty$.

Although the system \eqref{eqn:wna_adj_W_problem} appears overdetermined, Eq.~\eqref{eqn:wna_adj_bif_condition} turns out to be a restatement of \eqref{eqn:wna_eta_alg_eqn}.  To see this, we multiply \eqref{eqn:wna_adj_W_eqn} by $\eta(u)$, which satisfies \eqref{eqn:wna_Odelta_eqn}, and integrate over $u$. This yields
\begin{align}
	\alpha_0\int_0^\infty u \eta(u)\,\mathrm{d}u &= \int_0^\infty\left[\varepsilon\dft{\eta}{u} + \dfo{}{u}\left[\lambda(u - u^*)\eta\right] - \mathcal{D}_0k^2\eta \right]W\,\mathrm{d}u \nonumber \\
	&= -\int_0^\infty \fc \dfo{n^*}{u}W(u)\,\mathrm{d}u,
\end{align}
and we conclude that Eqs.~\eqref{eqn:wna_adj_bif_condition} and \eqref{eqn:wna_eta_alg_eqn} are equivalent. This implies that \eqref{eqn:wna_adj_W_problem} can only have a solution at the bifurcation point.  

Next we obtain solutions of \eqref{eqn:wna_adj_W_eqn} in the regime $\varepsilon\ll 1$. Following our analysis of \eqref{eqn:canonical_q_ODE} for $q(u)$, we expand $W$ as
\begin{equation}
	W(u) \sim \sum_{j=0}^\infty W_j(u)\varepsilon^j,
	\label{eqn:wna_W_asymp_expand}
\end{equation}
and substitute into \eqref{eqn:wna_adj_W_eqn} to obtain the following set of ODEs for the coefficients
\begin{subequations}
\begin{align}
	\lambda (u - u^*) W_0' + \mathcal{D}_0(u)k^2W_0  &=\alpha_0 u , \label{eqn:wna_w0_eqn}\\
	\lambda (u - u^*) W_j' + \mathcal{D}_0(u)k^2W_j  &= W_{j-1}'', \quad j=1,\ldots.
\end{align}
\label{eqn:wna_wj_eqns}%
\end{subequations}
Since the only difference between \eqref{eqn:canonical_q0_q1_ode} and \eqref{eqn:wna_wj_eqns} is the inhomogeneous term on the RHS of \eqref{eqn:wna_w0_eqn}, our analysis of the Stokes phenomenon for $q$ also applies for $W$ (see Appendix \ref{app:stokes}). We therefore deduce that solutions of \eqref{eqn:wna_wj_eqns} must be smooth at $u^*$ in order for $W$ to be polynomially bounded as $u\to0,\infty$. This allows us to seek formal power series solutions for $W_j$ in the form
\begin{equation}
	W_j(u) = \sum_{\ell=0}^\infty \updown{W}{j}{\ell} (u - u^*)^\ell.
	\label{eqn:wna_wj_series}
\end{equation}
Using standard power series methods, we calculate the coefficients as
\begin{subequations}
\begin{align}
		\updown{W}{0}{0} &=\frac{\alpha_0 u^*}{\mathcal{D}_{0*}k^2},\\
		\updown{W}{0}{1}&=  \frac{\alpha_0 - \updown{W}{0}{0}\mathcal{D}_{0*}'k^2}{\mathcal{D}_{0*}k^2 + \lambda}, \\
		\updown{W}{0}{\ell} &= -\frac{k^2}{\mathcal{D}_{0*}k^2 + \ell \lambda }\sum_{m=1}^\ell \updown{W}{0}{\ell-m}\frac{\mathcal{D}_{0*}^{(m)}}{m!}, \quad \ell = 2,\ldots, \\
		\updown{W}{j}{\ell} &= \frac{(\ell + 2)(\ell + 1)\updown{W}{j-1}{\ell + 2} -k^2\sum\limits_{m=1}^\ell \updown{W}{j}{\ell - m}\frac{\mathcal{D}_{0*}^{(m)} }{m!}}{\mathcal{D}_{0*}k^2 + \ell \lambda}, ~j=1,\ldots, ~ \ell = 0,\ldots.
 \end{align}
 \label{eqn:wna_w_j_ps_coeffs}%
\end{subequations}
We will only require a finite number of the coefficients $\updown{W}{j}{\ell}$ in the course of our analysis, all of which can be calculated \emph{via} the recurrence relations \eqref{eqn:wna_w_j_ps_coeffs}. We remark that even though the nontrivial solutions of the homogeneous adjoint problem $\mathcal{L}^*(\mbf{w})=\mbf{0}$ are not exponentially localized, only their local behaviour will feed into our final result for the amplitude $A(\tau)$. Thus a power series solution for $W_j$ is sufficient.

Now that the nullspace of $\mathcal{L}^*$ is known, we can obtain the solvability conditions. The system $\mathcal{L}(\mbf{v})=\mbf{g}:=(g_1(x,u),g_2(x))^T$ with the boundary conditions \eqref{eqn:v_BCs} has a solution if, and only if, the two solvability conditions 
\begin{subequations}
    \begin{align}
  \int_{\Omega\times\mathbb{R}_+} g_1(x,u)\,\mathrm{d}x\mathrm{d}u &= \est,
	\label{eqn:wna_solv_cond1_gen} \\
	\int_{\Omega\times\mathbb{R}_+} g_1(x,u)W(u)\cos(kx)\,\mathrm{d}u\mathrm{d}x + \int_{\Omega} g_2(x)\cos(kx)\,\mathrm{d}x &= 0. 
    \label{eqn:wna_solv_cond2_gen} 
    \end{align}
    \label{eqn:wna_solv_cond_gen}%
\end{subequations}
are satisfied. The exponentially small terms in \eqref{eqn:wna_solv_cond1_gen} are due to the fact that $w_1\equiv 1$, $w_2\equiv 0$ only satisfies \eqref{eqn:WNA_adj_op}--\eqref{eqn:WNA_V_adj_BC} up to exponentially small corrections.

We demand that there exists a solution of the $\mathcal{O}(\delta^2)$ problem in \eqref{eqn:gov_canonical_wna_Odelta2}  and therefore impose the solvability conditions \eqref{eqn:wna_solv_cond_gen}. Substituting the RHS of \eqref{eqn:gov_canonical_wna_eta2_eqn} into \eqref{eqn:wna_solv_cond_gen}, we find that both conditions are trivially satisfied and no information on $A(\tau)$ is obtained, as expected. The governing equation for $A(\tau)$  will be deduced from the solvability conditions at  $\mathcal{O}(\delta^3)$. Since the solutions $\eta_2$ and $c_2$ will feed into the terms at order $\mathcal{O}(\delta^3)$, our next task is to determine the particular integral of  \eqref{eqn:gov_canonical_wna_eta2_eqn}.

As $\eta_1$ and $c_1$ are proportional to $A(\tau)\cos(kx)$, the RHS of \eqref{eqn:gov_canonical_wna_eta2_eqn} is proportional to $A^2\cos^2(kx)$. We therefore seek a particular integral of the form
\begin{equation}
	\begin{pmatrix}
	\eta_2(x,u,\tau) \\
    c_2(x,\tau)  
    \end{pmatrix}
    =  A^2(\tau)
    \begin{pmatrix}
       \eta_{20}(u) + \eta_{22}(u)\cos(2kx) \\
        c_{20} + c_{22}\cos(2kx),
    \end{pmatrix}
	\label{eqn:wna_eta2_c2_mode}
\end{equation}
where $\eta_{ij}$, $c_{ij}$ are new unknowns to be determined. We omit contributions from the homogeneous solution $\cos(kx)(\eta(u),1)^T$, as these terms will not affect the solvability condition at next order. 

We obtain governing equations for the unknowns $\eta_{ij}$, $c_{ij}$ by substituting \eqref{eqn:wna_eta2_c2_mode} into \eqref{eqn:gov_canonical_wna_eta2_eqn} and using the linear independence of the Fourier modes. For $\eta_{20}$ and $c_{20}$ we obtain
\begin{subequations}
\begin{align}
	\varepsilon \dft{\eta_{20}}{u} - \dfo{}{u}\left[f(u,c^*)\eta_{20}\right] - c_{20}\fc\dfo{n^*}{u} &= \frac{1}{2}\left[\frac{\fcc}{2}\dfo{n^*}{u} + \fc\dfo{\eta}{u}\right], \label{eqn:wna_Odelta2_eta20_eqn}\\
	\beta c_{20} - \alpha_0\int_0^\infty u \eta_{20}(u)\,\mathrm{d}u &=0, \label{eqn:wna_Odelta2_c20_eqn}
\end{align}
\label{eqn:wna_Odelta2_cos0_prob}%
\end{subequations}
and similarly for $\eta_{22}$ and $c_{22}$:
\begin{subequations}
\begin{align}
	\varepsilon \dft{\eta_{22}}{u} - \dfo{}{u}\left[f(u,c^*)\eta_{22}\right] - c_{22}\fc\dfo{n^*}{u} - 4\mathcal{D}_0k^2\eta_{22} &= \frac{1}{2}\left[\frac{\fcc}{2}\dfo{n^*}{u} + \fc\dfo{\eta}{u}\right], \label{eqn:wna_Odelta2_eta22_eqn}\\
	\left(4D_ck^2 + \beta\right)c_{22} - \alpha_0\int_0^\infty u \eta_{22}(u)\,\mathrm{d}u & = 0. \label{eqn:wna_Odelta2_c22_eqn} 
\end{align}
\label{eqn:wna_Odelta2_cos2kx_prob}%
\end{subequations}
Both $\eta_{20}$ and $\eta_{22}$ are subject to no-flux boundary conditions in $u$ (cf.~\eqref{eqn:gov_canonical_eta2_BC}--\eqref{eqn:gov_canonical_c2_BC}). 

Eq.~\eqref{eqn:wna_Odelta2_eta20_eqn} has the following exact solution for $\eta_{20}$:
\begin{subequations}
\begin{align}
    \eta_{20}(u) &= r(u)\exp\left(-\frac{\lambda (u-u^*)^2}{2\varepsilon}\right),  \label{eqn:wna_r_eqn} \\
    r(u) &=  \left(c_{20}\fc + \frac{\fcc}{4}\right)\rho^*\sqrt{\frac{\lambda}{2\pi\varepsilon^3}} (u - u^*) + \varepsilon^{-5/2}\frac{\fc}{2} Q(u) + \bar{r}\sqrt{\frac{\lambda}{2\pi\varepsilon}}, 
	\label{eqn:wna_r_sol} \\
    Q(u)&:= \int_{u^*}^u q(u)\,\mathrm{d}u.
\end{align}
Here $\bar{r}$ is an arbitrary integration constant that we fix by imposing the normalisation condition \eqref{eqn:normalisation}, i.e.~the perturbation must have zero mass: $\int_{\Omega\times\mathbb{R}_+}\eta_2\,\mathrm{d}u\mathrm{d}x=0$.  Thus
\begin{equation}
    \bar{r}= -\frac{\fc}{2\varepsilon^{5/2}}\int_0^\infty Q(u) \exp\left(-\frac{\lambda(u - u^*)^2}{2\varepsilon}\right) + \est\label{eqn:wna_rbar_exact}
\end{equation}
\label{eqn:wna_eta_20_soln}%
\end{subequations}
We apply Laplace's method to evaluate the integrals in \eqref{eqn:wna_Odelta2_c20_eqn} and \eqref{eqn:wna_rbar_exact} and find that $c_{20}$ is given by
\begin{equation}
c_{20} = \frac{\alpha_0\left(\rho^*\fcc + \sqrt{\frac{8\pi}{\lambda^3}}\fc\left(\lambda \updown{q}{1}{0} + \updown{q}{0}{2}\right)\right)}{4\left(\lambda\beta - \alpha_0\rho^*\fc\right)} + \mathcal{O}(\varepsilon).
\label{eqn:wna_c20_val}
\end{equation}
In principle $c_{20}$ can be expressed in terms of the original model parameters by eliminating $\updown{q}{1}{0}$ and  $\updown{q}{2}{0}$ with \eqref{eqn:canonical_wna_qj_coeffs}.

The next step in solving \eqref{eqn:gov_canonical_wna_eta2_eqn} is to determine  $\eta_{22}$ and $c_{22}$ from \eqref{eqn:wna_Odelta2_cos2kx_prob}. As before, we seek a WKBJ solution in the form
\begin{equation}
	\eta_{22}(u) = \varepsilon^{-5/2}s(u)\exp\left(-\frac{\lambda(u - u^*)^2}{2\varepsilon}\right),
	\label{eqn:wna_eta22_wkb}
\end{equation}
where we include the $\varepsilon^{-5/2}$ scaling with hindsight so that $s=\mathcal{O}(1)$ as $\varepsilon\to 0$, as seen from \eqref{eqn:wna_s_eqn} below. From \eqref{eqn:wna_Odelta2_eta22_eqn}, the amplitude equation for $s$ is then
\begin{equation}
	\varepsilon s'' - \lambda (u - u^*) s' - 4\mathcal{D}_0k^2 s =  \frac{\fc}{2}\left(\varepsilon q' - \lambda (u - u^*)q\right) -\varepsilon\rho^*\sqrt{\frac{\lambda^3}{2\pi}}\left(\fc c_{22} + \frac{\fcc}{4}\right).
	\label{eqn:wna_s_eqn}
\end{equation}
Since $q=\mathcal{O}(1)$, as seen from \eqref{eqn:canonical_q_expand}, \eqref{eqn:canonical_qj_ps}--\eqref{eqn:wna_q0_ps_coeffs}, the RHS of \eqref{eqn:wna_s_eqn} is also $\mathcal{O}(1)$, so we expand $s$ in a regular asymptotic series as
\begin{equation}
	s(u)\sim \sum_{j=0}^\infty s_j(u)\varepsilon^j.
	\label{eqn:wna_s_asymp_expand}
\end{equation}
Inserting \eqref{eqn:wna_s_asymp_expand} into \eqref{eqn:wna_s_eqn} and collecting powers of $\varepsilon$, we obtain the following set of ODEs for $s_j$
\begin{subequations}
	\begin{align}
		\lambda (u - u^*)s_0' + 4\mathcal{D}_0k^2 s_0 =&\frac{\fc}{2}\lambda (u - u^*)q_0(u), \\
		\lambda (u - u^*)s_j' + 4\mathcal{D}_0k^2 s_j =& s_{j-1}'' - \frac{\fc}{2}\left(q_{j-1}' - \lambda (u -u^*)q_j\right) \nonumber \\
		& + \delta_{j1}\rho^*\sqrt{\frac{\lambda^3}{2\pi}}\left(\fc c_{22} + \frac{\fcc}{4}\right)(u - u^*), \quad j=1,\ldots
	\end{align}
	\label{eqn:wna_sj_eqns}%
\end{subequations}
As before, our Stokes phenomenon argument implies that $s$ is smooth at $u=u^*$ (see Appendix \ref{app:stokes}). Hence, the solution satisfying the boundary conditions can be written as a formal power series
\begin{equation}
	s_j(u)=\sum_{\ell=0}^\infty \updown{s}{j}{\ell}(u - u^*)^\ell, \qquad j=0,\ldots.
	\label{eqn:wna_sj_series}
\end{equation}
In a similar way to before, we calculate the coefficients with standard power series methods and find, after some algebra,
\begin{subequations}
\begin{align}
	\updown{s}{0}{0} =& \updown{s}{0}{1} = 0, \\
	\updown{s}{0}{\ell} =& \frac{\frac{1}{2}\fc\updown{q}{0}{\ell - 1} - 4k^2\sum_{m=1}^\ell \frac{\mathcal{D}_{0*}^{(m)}}{m!}\updown{s}{0}{\ell - m}}{4\mathcal{D}_{0*}k^2 + \ell \lambda}, \quad \ell = 2,\ldots, \\
	\updown{s}{j}{0} =& \frac{2\updown{s}{j-1}{2}}{4\mathcal{D}_{0*}k^2}, \quad j=1,\ldots, \\
	\updown{s}{j}{\ell} =& \bigg[(\ell + 1)(\ell + 2)\updown{s}{j-1}{\ell} - 4k^2\sum_{m=1}^\ell\frac{\mathcal{D}_{0*}}{m!}\updown{s}{j}{\ell-m} - \frac{\fc}{2}\left[(\ell + 1)\updown{q}{j-1}{\ell+1} - \lambda \updown{q}{j}{\ell - 1}\right] \nonumber \\
	&+ \delta_{\ell1}\delta_{j1}\rho^*\sqrt{\frac{\lambda^3}{2\pi}}\left(\fc c_{22} + \frac{\fcc}{4}\right)\bigg]\left(4\mathcal{D}_{0*}k^2 + \ell \lambda\right)^{-1}, \quad \ell,j=1,\ldots. \label{eqn:wna_s_ps_sjl}
\end{align}
\label{eqn:wna_s_ps_coeffs}%
\end{subequations}
In calculating the coefficients in \eqref{eqn:wna_s_ps_coeffs}, we have used the expansions for $q_j$ in \eqref{eqn:canonical_qj_ps} to write the coefficients in terms of those for $q_j$. 

The last step at $\mathcal{O}(\delta^2)$ is to determine $c_{22}$ from \eqref{eqn:wna_Odelta2_c22_eqn}. Since $\eta_{22}$ explicitly depends on $c_{22}$ through the coefficients of the amplitude $s$ (see \eqref{eqn:wna_s_ps_sjl}), we simplify the calculation by eliminating $\eta_{22}$ from the integrand in \eqref{eqn:wna_Odelta2_c22_eqn}. We accomplish this by observing that the adjoint of \eqref{eqn:wna_Odelta2_eta22_eqn} is \eqref{eqn:wna_adj_W_eqn}, only with the replacement $k\mapsto 2k$. We therefore have
\begin{equation}
	\alpha_0 \int_0^\infty u\eta_{22}\,\mathrm{d}u
	=\int_0^\infty\left[\left(c_{22} \fc + \frac{\fcc}{4}\right)n^* + \frac{\fc}{2}\eta\right]\tilde{W}\,\mathrm{d}u,
	\label{eqn:wna_adj_int_relation}
\end{equation}
where $\tilde{W}$ denotes $W$ with the mapping $k\mapsto 2k$. Evaluating the integral using \eqref{eqn:n_ss} for $n^*$; \eqref{eqn:wna_Odelta_wkbj}, \eqref{eqn:canonical_q_expand}, and \eqref{eqn:canonical_wna_qj_coeffs} for $\eta$;  \eqref{eqn:wna_W_asymp_expand} and \eqref{eqn:wna_wj_series} for $\tilde{W}$; and Laplace's method, we find
\begin{equation}
	c_{22} = \frac{\rho^*\fcc\updown{\tilde{W}}{0}{1} + 2\fc\sqrt{\frac{2\pi}{\lambda^3}}\left[\left(\lambda\updown{q}{1}{0} + \updown{q}{0}{2}\right)\updown{\tilde{W}}{0}{1} + + 2\updown{q}{0}{1}\updown{\tilde{W}}{0}{2}\right]}{4\left(4D_ck^2 +\beta  - \rho^*\fc \updown{\tilde{W}}{0}{1}\right)} + \mathcal{O}(\varepsilon),
	\label{eqn:wna_c22_val}
\end{equation}
where $\updown{\tilde{W}}{j}{\ell}$ denotes $\updown{W}{j}{\ell}$ (given in \eqref{eqn:wna_w_j_ps_coeffs}) with $k\mapsto 2k$.

We have now determined $\eta_{20}$, $c_{20}$, $\eta_{22}$, and $c_{22}$ in terms of $A(\tau)$ (given in \linebreak Eqs.~\eqref{eqn:wna_eta_20_soln}--\eqref{eqn:wna_eta22_wkb}, \eqref{eqn:wna_s_asymp_expand}, and \eqref{eqn:wna_sj_series}--\eqref{eqn:wna_c22_val}). Thus $\eta_2$ and $c_2$ in \eqref{eqn:wna_eta2_c2_mode} are known in terms of $A(\tau)$, which completes our analysis at $\mathcal{O}(\delta^2)$. Our final task is to perform an analysis at $\mathcal{O}(\delta^3)$ which will yield the amplitude equation for $A(\tau)$ that we seek.

\subsection{Analysis at  $\mathcal{O}(\delta^3)$} \label{sec:Odelta3}

Finally, we investigate the governing equations at $\mathcal{O}(\delta^3)$. Inserting \eqref{eqn:wna_expand} and the scalings in \eqref{eqn:WNA_Dp_perturb}, \eqref{eqn:WNA_t_scaling} into the governing equations \eqref{eqn:gov_eqs}--\eqref{eqn:gov_BCs} and collecting  terms at $\mathcal{O}(\delta^3)$ yields
\begin{subequations}
\begin{equation}
	\mathcal{L}\begin{pmatrix}
		\eta_3\\
		c_3
	\end{pmatrix} 
	=
	\begin{pmatrix}
	\left(c_1c_2\fcc + c_1^3 \frac{\fccc}{6}\right)\dfo{n^*}{u} + \left(c_2\fc + c_1^2\frac{\fcc}{2}\right)\pfo{\eta_1}{u}  
	 + c_1\fc\pfo{\eta_2}{u} + \pfo{\eta_1}{\tau} -d\pft{\eta_1}{x} \\
	 \pfo{c_1}{\tau}
	\end{pmatrix}.
    \label{eqn:gov_eta3_c3}
\end{equation} 
\begin{align}
	\pfo{c_3}{x}&=\pfo{\eta_3}{x} =0, \quad x=0,L,\\
	\varepsilon\pfo{\eta_3}{u} - f(u,c^*)\eta_3 &\to 0, \quad \text{ as } u\to 0,\infty.
\end{align}
\label{eqn:gov_delta3}
\end{subequations}

To obtain the amplitude $A(\tau)$, we do not need to solve \eqref{eqn:gov_delta3}; all we require is that a solution exists. This means we can impose the two solvability conditions \eqref{eqn:wna_solv_cond_gen}, which will yield an amplitude equation for $A(\tau)$.  The condition \eqref{eqn:wna_solv_cond1_gen} says that integral of the first component of the RHS of  \eqref{eqn:gov_eta3_c3} is exponentially small, which is trivially satisfied due to the orthogonality of the Fourier modes (we recall that $c_1$ and $\eta_1$ are proportional to $\cos(kx)$ while $c_2$ are $\eta_2$ are a linear combination of a constant and $\cos(2kx)$). As flagged previously, even if we had included a contribution from the homogeneous solution $\cos(kx)(\eta(u),1)^T$ to $(\eta_2,c_2)^T$, the integral would still be exponentially small,  and the solvability condition \eqref{eqn:wna_solv_cond1_gen} would still be satisfied.

Imposing the second solvability condition \eqref{eqn:wna_solv_cond2_gen} will yield the weakly non-linear form for the bifurcation. The integrals over $x$ in \eqref{eqn:wna_solv_cond2_gen} can be evaluated using   \eqref{eqn:wna_Odelta_separation} for $\eta_1$ and $c_1$, \eqref{eqn:wna_adj_mode} for $w_1$ and $w_2$, as well as \eqref{eqn:wna_eta2_c2_mode} for $\eta_2$ and $c_2$. This gives
\begin{subequations}
\begin{align}
	\dfo{A}{\tau} &= - \frac{I_1}{I_0 + 1} k^2 A - \frac{\mu}{I_0 + 1} A^3, \label{eqn:wna_C_ODE} \\
    \mu&:=\left(c_{20} + \frac{c_{22}}{2}\right)\left(\fcc I_2 + \fc I_3 \right) + \frac{\fccc I_2 + 3\fcc I_3}{8} + \fc\left(I_4 + \frac{I_5}{2}\right)
    \label{eqn:wna_mu_def}
\end{align}
where the integrals $I_0,\ldots, I_5$ are given by
\begin{IEEEeqnarray}{rll}
	I_0&:=\int_0^\infty \eta(u) W(u)\,\mathrm{d}u, \quad 
	&I_1:=\int_0^\infty d(u)\eta(u) W(u)\,\mathrm{d}u, \\
	I_2&:=\int_0^\infty \dfo{n^*}{u}W\,\mathrm{d}u, \quad
	&I_3:=\int_0^\infty \dfo{\eta}{u} W\,\mathrm{d}u,\\ 
	I_4&:=\int_0^\infty \dfo{\eta_{20}}{u} W\,\mathrm{d}u, \quad
	&I_5:=\int_0^\infty \dfo{\eta_{22}}{u} W \,\mathrm{d}u. 
\end{IEEEeqnarray}
\label{eqn:wna_Ceqn_mu_I_de}%
\end{subequations}
The leading-order contributions to the integrals $I_j$ in \eqref{eqn:wna_Ceqn_mu_I_de} can be explicitly evaluated with Laplace's method  since $n^*$, $\eta$, $\eta_{20}$, and $\eta_{22}$ are exponentially localized. A lengthy but straightforward calculation gives
\begin{subequations}
\begin{align}
	I_0=& \frac{\alpha_0\rho^*\fc\left[\mathcal{D}_{0*}k^2 - \left(2\mathcal{D}_{0*}k^2 + \lambda\right)u^*\frac{\mathcal{D}_{0*}'}{\mathcal{D}_{0*}}\right]}{\left(\mathcal{D}_{0*}k^2 + \lambda\right)^2\mathcal{D}_{0*}k^2} + \mathcal{O}(\varepsilon), \label{eqn:wna_I0}\\
	I_1=& \frac{\alpha_0\rho^* u^* \fc}{\mathcal{D}_{0*}k^2\left(\mathcal{D}_{0*}k^2 + \lambda\right)}d'_* + \mathcal{O}(\varepsilon), \label{eqn:wna_I1}\\
	I_2=& \frac{\alpha_0\rho^*\left(u^*\frac{\mathcal{D}_{0*}'}{\mathcal{D}_{0*}} - 1\right)}{\mathcal{D}_{0*}k^2 + \lambda} + \mathcal{O}(\varepsilon), \\\
	I_3=& - \sqrt{\frac{2\pi}{\lambda}}\left[\left(\updown{q}{1}{0} + \frac{\updown{q}{0}{2}}{\lambda}\right)\updown{W}{0}{1} + \frac{2}{\lambda}\updown{q}{0}{1}\updown{W}{0}{2} \right] + \mathcal{O}(\varepsilon),\\
	I_4=& -\sqrt{\frac{2\pi}{\lambda^5}}\fc\left[\left(\updown{q}{0}{2} + \lambda \updown{q}{1}{0}\right)\updown{W}{0}{2} + \frac{3}{2}\updown{q}{0}{1}\updown{W}{0}{3}\right] \nonumber \\ &- \frac{\rho^*}{\lambda}\left(2c_{20} \fc + \frac{\fcc}{2}\right)\updown{W}{0}{2} + \mathcal{O}(\varepsilon),\\
	I_5=& -\sqrt{\frac{2\pi}{\lambda}}\Bigg[\left(\updown{s}{2}{0} + \frac{\updown{s}{1}{2}}{\lambda} + \frac{3\updown{s}{0}{4}}{\lambda^2}\right)\updown{W}{0}{1} + \frac{2}{\lambda}\left(\updown{s}{1}{1} + \frac{3\updown{s}{0}{3}}{\lambda}\right)\updown{W}{0}{2} \nonumber \\
	&+ \frac{3}{\lambda}\left(\updown{s}{1}{0} + \frac{3\updown{s}{0}{2}}{\lambda}\right)\updown{W}{0}{3}\Bigg] + \mathcal{O}(\varepsilon),
\end{align}
\label{eqn:wna_integrals_simp}%
\end{subequations}
where we have used Eq.~\eqref{eqn:canonical_density_ss}--\eqref{eqn:canonical_density_ss_normalisation} for $n^*$, Eqs.~\eqref{eqn:wna_Odelta_wkbj}, \eqref{eqn:canonical_q_expand}, and \eqref{eqn:canonical_qj_ps} for $\eta$, Eqs.~\eqref{eqn:wna_W_asymp_expand}, and \eqref{eqn:wna_wj_series} for $W$, Eq.~\eqref{eqn:wna_eta_20_soln} for $\eta_{20}$, and finally \eqref{eqn:wna_eta22_wkb}, \eqref{eqn:wna_s_asymp_expand}, and \eqref{eqn:wna_sj_series} for $\eta_{22}$. We note that if $d_*:=d(u^*)$ were non-zero (see \eqref{eqn:d_perturb}), then there would be an additional term in $I_1$ proportional to $d_*$. Since both $d_*$ and $d_*'$ vanish at the bifurcation point defined by $\theta_0=0$ (see \eqref{eqn:theta_perturb}--\eqref{eqn:WNA_Dp_perturb}), the assumption \eqref{eqn:d_perturb}  does not fundamentally alter the bifurcation structure. We also note that the integral $I_0$ is positive because $\mathcal{D}_*'<0$, so the denominators in \eqref{eqn:wna_C_ODE} are nonzero. The integrals $I_j$ can, in principle, be expressed in terms of the original model parameters with the recursion relations \eqref{eqn:canonical_wna_qj_coeffs}, \eqref{eqn:wna_w_j_ps_coeffs}, and \eqref{eqn:wna_s_ps_coeffs}, together with the results  \eqref{eqn:wna_c20_val} and \eqref{eqn:wna_c22_val} for $c_{20}$ and $c_{22}$. For the sake of brevity, we simplify only the linear terms in  Eq.~\eqref{eqn:wna_C_ODE} using \eqref{eqn:wna_I0} and \eqref{eqn:wna_I1}, leading to
\begin{equation}
	\dfo{A}{\tau} \sim  \frac{-\alpha_0 \rho^* \fc u^*\left(\mathcal{D}_{0*}k^2 + \lambda\right)d_*'k^2}{\mathcal{D}_{0*}k^2\left(\mathcal{D}_{0*}k^2 + \lambda\right)^2 + \alpha_0 \rho^*\fc\left[ \mathcal{D}_{0*}k^2 - \left(2\mathcal{D}_{0*}k^2 + \lambda\right)u^*\frac{\mathcal{D}_{0*}'}{\mathcal{D}_{0*}}\right]}A -\frac{\mu}{I_0+1}A^3,
	\label{eqn:wna_C_ODE_simp}
\end{equation}
where we neglect $\mathcal{O}(\varepsilon)$ terms in the coefficients. We recall that the amplitude $A(\tau)$ governs the slow-time dependence of the perturbations $\eta_1(x,u,\tau)$ and $c_1(x,\tau)$ in \eqref{eqn:wna_expand}, \eqref{eqn:wna_Odelta_separation}. We remark that we have used \eqref{eqn:d_perturb}, i.e.~$d(u^*)=0$ in our evaluation of $I_1$.

Eq.~\eqref{eqn:wna_C_ODE_simp} governs the weakly non-linear behaviour of the system near $d_*'=0$ and describes a pitchfork bifurcation.   To understand the local behaviour of the non-uniform steady states in terms of the cell density $\rho$, we eliminate the local variables $d_*'$ and $A$ in \eqref{eqn:wna_C_ODE_simp} in favour of the global bifurcation parameter $\mathcal{D}_*'$ and cell density $\rho$ using Eqs.~\eqref{eqn:rho_def}, \eqref{eqn:WNA_Dp_perturb}, and \eqref{eqn:wna_expand}. Thus as  $\mathcal{D}_*' \to \mathcal{D}_{0*}'$, we find that the non-uniform steady states of \eqref{eqn:wna_C_ODE_simp} have the asymptotic behaviour
\begin{equation}
	\rho(x) =\rho^* \pm \mathcal{D}_{0*}'\sqrt{-\frac{\alpha_0u^*(\rho^*\fc)^3\left(\mathcal{D}'_* - \mathcal{D}_{0*}'\right)}{\mu \mathcal{D}_{0*}^3\left(\mathcal{D}_{0*}k^2 + \lambda\right)^3}} \cos(kx)\big(1 + \mathcal{O}(\varepsilon)\big) + \mathcal{O}(\mathcal{D}_*'-\mathcal{D}_{0*}'),
	\label{eqn:wna_rho_nonuniform}
\end{equation}
whenever $\mu$ and $\mathcal{D}'_* - \mathcal{D}_{0*}'$ have opposite signs. The sign of $\mu$ determines the linear stability of the non-uniform steady state \eqref{eqn:wna_rho_nonuniform}. We recall that the uniform steady state $\rho=\rho^*$ is stable for $\mathcal{D}'_* >\mathcal{D}_{0*}'$ and unstable for $\mathcal{D}'_* <\mathcal{D}_{0*}'$. Thus the pitchfork bifurcation is subcritical when $\mu<0$ and supercritical when $\mu>0$. The condition 
\begin{equation}
    \mu:=\left(c_{20} + \frac{c_{22}}{2}\right)\left(\fcc I_2 + \fc I_3 \right) + \frac{\fccc I_2 + 3\fcc I_3}{8} + \fc\left(I_4 + \frac{I_5}{2}\right)=0,
	\label{eqn:wna_sub_super_cond}
\end{equation}
therefore defines the transition between subcritical and supercritical pitchfork bifurcations. Eqs.~\eqref{eqn:wna_C_ODE_simp}--\eqref{eqn:wna_sub_super_cond} are the main results of our weakly nonlinear analysis.

\section{Numerical Results} \label{sec:numerical_results}

To demonstrate the predictions of our weakly nonlinear theory, we generate numerical solutions of the governing equations \eqref{eqn:gov_eqs}--\eqref{eqn:rho_def} at steady state. To facilitate these calculations, we make the following convenient choices for the diffusion coefficient $\mathcal{D}(u)$ and  reaction kinetics $g(c)$ 
\begin{subequations}
    \begin{align}
        \mathcal{D}(u) &= \mathcal{D}_*  - (\mathcal{D}_* - \mathcal{D}_\infty)\tanh w(u), \quad w(u)=-\frac{\mathcal{D}_*'}{\mathcal{D}_* - \mathcal{D}_\infty}(u - u^*), \label{eqn:D_eg}\\
        g(c) &= a + \frac{Vc}{K+c}.
    \end{align}
    \label{eqn:gD_specific}%
\end{subequations}
Here $\mathcal{D}_*$ and $\mathcal{D}_\infty$ are fixed parameters, while $\mathcal{D}_*'$ is our bifurcation parameter (a shifted version of $\theta$ used in our analysis). The form \eqref{eqn:D_eg} for $\mathcal{D}$ satisfies the assumption \eqref{eqn:d_perturb}.  Our choice for $g(c)$ in \eqref{eqn:gD_specific} is the same as the form used in \cite{ridgway2023motility}, which incorporates positive feedback -- a canonical feature of quorum sensing systems \cite{papenfort2016quorum,dalwadi2021emergent}. We use the open source library oomph-lib \cite{heil2006oomph} for our numerical calculations (see \cite{ridgway2023motility} for details and codes). Briefly, we use a  Galerkin finite element method with quadratic Lagrange elements to compute numerical solutions, along with pseudo-arclength continuation to track solutions as we vary the bifurcation parameter $\mathcal{D}_*'$.

\begin{figure}
    \centering
    \includegraphics[width=0.99\textwidth]{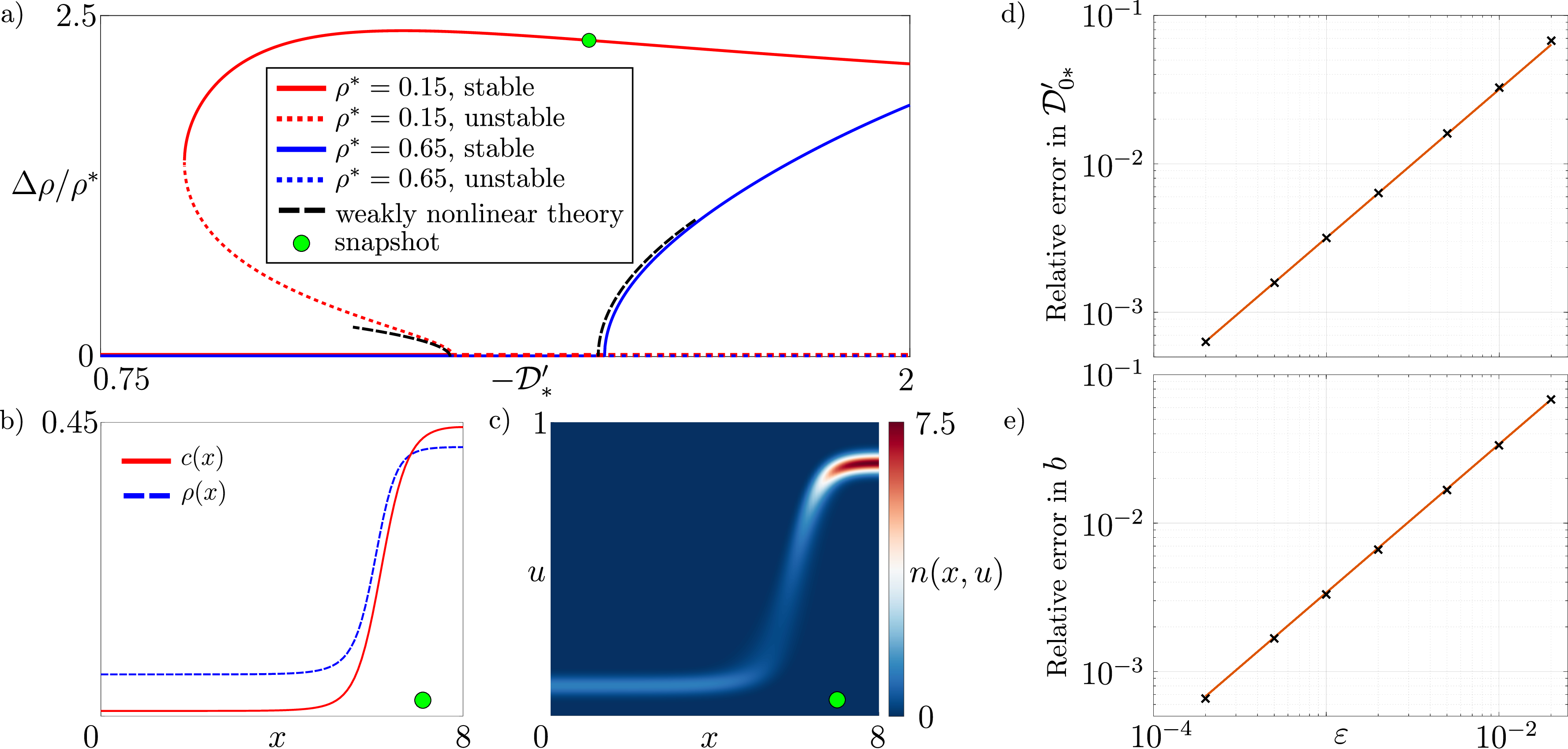}
    \caption{a): sample bifurcation diagrams when the pitchfork bifurcation is supercritical and subcritical. Local nontrivial branches (black curves) from weakly nonlinear theory are shown near the pitchfork bifurcation. Red and blue curves with $\Delta\rho>0$ are computed from numerical solutions; the rest is from theory. Steady states are computed \emph{via} continuation, while stability is inferred from time-dependent calculations. b) and c): representative steady state profile for $\mathcal{D}_*'\approx -1.5$, and $\varepsilon=0.002$ (green circle, top left panel), computed numerically. Systematic comparison of theory with numerical results for $\rho^*=0.65$ and $\varepsilon\to0^+$. Relative errors in the leading order term of $\mathcal{D}_{0*}'$ in \eqref{eqn:bif_condition} d) and the coefficient $b$ in \eqref{eqn:b_def} e).  We plot a straight line passing through the origin to demonstrate that the error is $\mathcal{O}(\varepsilon)$. The numerical value of $b$ is determined by fitting a square root to the numerically computed points on the nontrivial branch near the pitchfork bifurcation.}
    \label{fig:bif_diagram_and_comparison}
\end{figure}

First, we construct bifurcation diagrams from our numerical solutions and compare the results with our theory. To this end we define $\Delta\rho:=|\rho(L)-\rho(0)|$ as an order parameter quantifying the degree of phase separation. This corresponds to the difference between the maximum and minimum cell densities near the pitchfork bifurcation. From the theoretical branch Eq.~\eqref{eqn:wna_rho_nonuniform}, we calculate
\begin{equation}
    \frac{\Delta\rho}{\rho^*} = b\sqrt{\left|\mathcal{D}_*'-\mathcal{D}_{0*}'\right|} + \mathcal{O}(\mathcal{D}_*'-\mathcal{D}_{0*}), \quad b=-2\mathcal{D}_{0*}'\sqrt{\frac{\alpha_0\rho^*u^*\fc^3}{|\mu| \mathcal{D}_{0*}^3\left(\mathcal{D}_{0*}k^2 + \lambda\right)^3}} + \mathcal{O}(\varepsilon),
    \label{eqn:b_def}
\end{equation}
whenever $\mu$ and $\mathcal{D}_*' - \mathcal{D}_{0*}'$ have opposite signs, where $\mu$ is defined in \eqref{eqn:wna_mu_def}. The nontrivial solution branches predicted by our weakly nonlinear theory agree well with the numerically computed bifurcation diagrams, as we show in Figure \ref{fig:bif_diagram_and_comparison} for both the supercritical and subcritical cases. Our weakly nonlinear theory predicts nontrivial solutions branches near the bifurcation, which are linearly stable when the bifurcation is supercritical, and linearly unstable when the bifurcation is subcritical.  In the supercritical case, the system admits small amplitude patterns near the  bifurcation when $\mathcal{D}_*'<\mathcal{D}_{0*}'$, which are well-approximated by our weakly nonlinear theory. In the subcritical case, the linearly stable steady state patterns have a large amplitude and are in general not accessible in a weakly nonlinear analysis. Our numerical calculations for the subcritical case show that the system relaxes to a large amplitude pattern, as seen in panels b) and c) of Figure \ref{fig:bif_diagram_and_comparison}.

\begin{figure}
    \centering
    \includegraphics[width=0.75\textwidth]{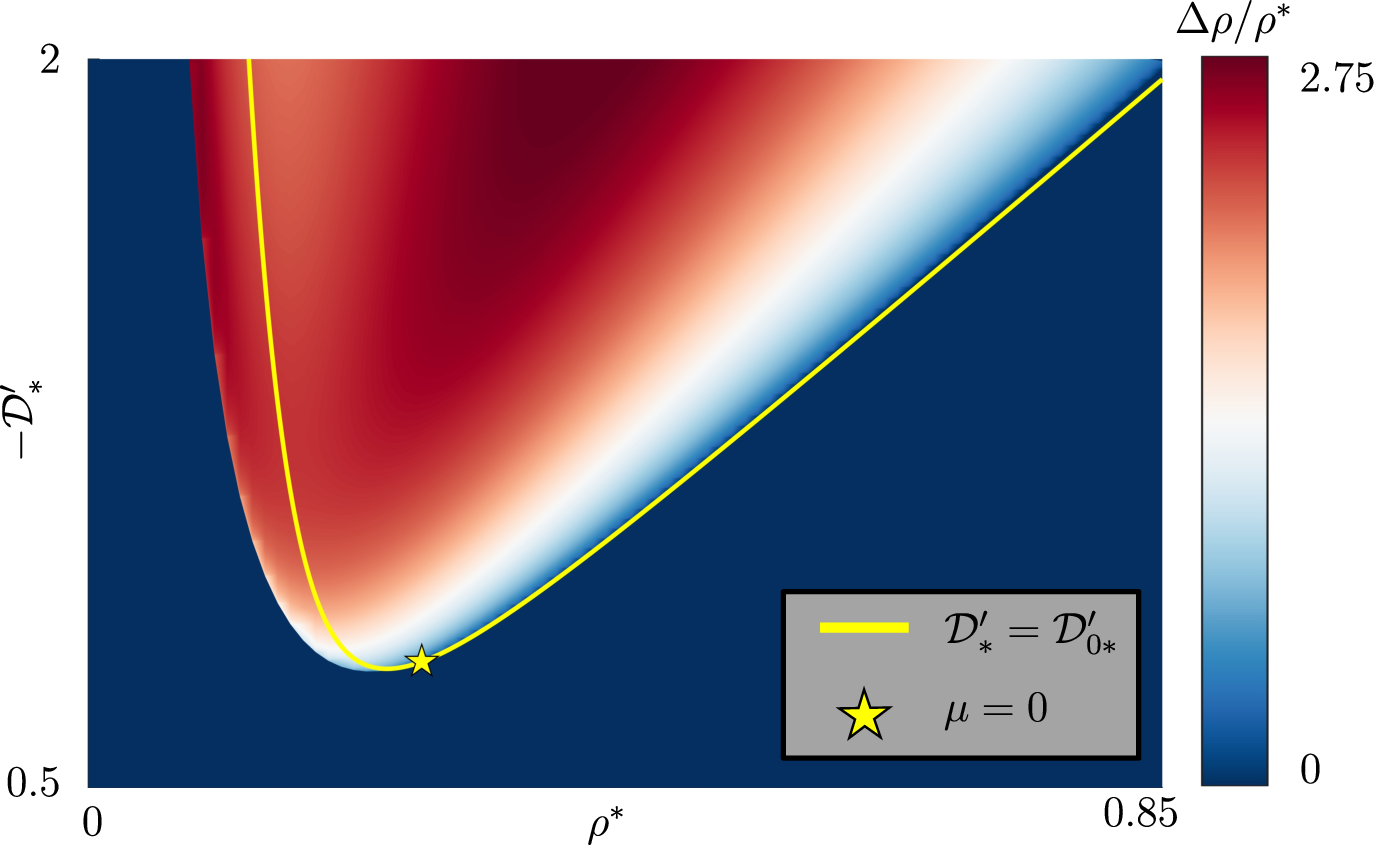}
    \caption{Phase diagram showing the order parameter $\Delta\rho/\rho^*$ in the phase-separated state ($\Delta\rho=0$ for the uniform solution). The yellow curve (Eq.~\eqref{eqn:bif_condition1}) traces out the bifurcation point; the spatially uniform state is linearly unstable above this curve and stable below. The bifurcation becomes subcritical as $\rho^*$ is decreased through $\rho^*\approx0.265$, defined by $\mu=0$  in \eqref{eqn:wna_sub_super_cond} (yellow star). The area wherein $\Delta\rho>0$ below the yellow curve denotes the region where phase separation is dynamically accessible (\emph{via} finite amplitude perturbation), but the uniform state is linearly stable.   }
    \label{fig:phase_diagram}
\end{figure}

The numerical bifurcation point and the local nontrivial branches agree reasonably well with the leading order terms in \eqref{eqn:bif_condition} and \eqref{eqn:b_def}, as shown in Figure 1. To verify our theory more systematically, we compare the numerically computed bifurcation point and coefficient $b$ with the theoretical results in \eqref{eqn:bif_condition} and \eqref{eqn:b_def} in the limit $\varepsilon\to0^+$.  The errors in both appear to be  $\mathcal{O}(\varepsilon)$ as $\varepsilon\to0^+$, as shown in Figure \ref{fig:bif_diagram_and_comparison}.

\begin{figure}
    \centering
    \includegraphics[width=0.95\linewidth]{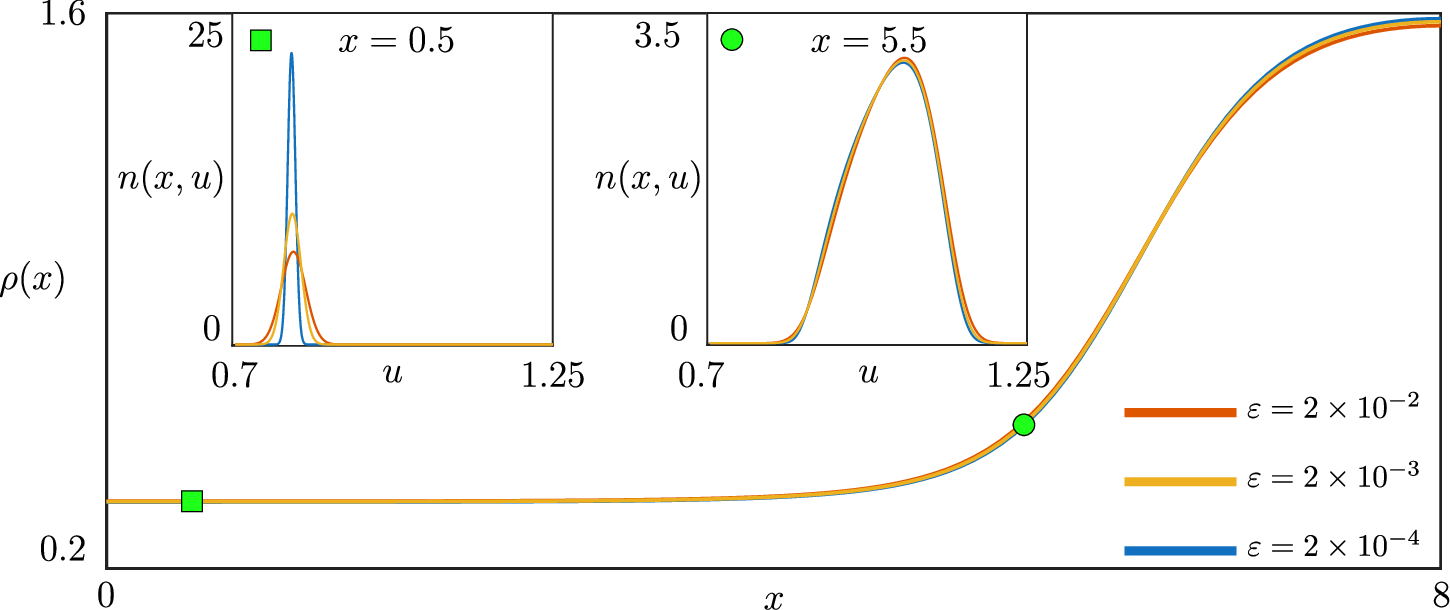}
    \caption{Numerically computed, non-uniform, steady state profile $\rho(x)$ and slices of $n(x,u)$ at fixed $x$, $(\rho^*,\mathcal{D}_*')=(0.65,-2)$. As $\varepsilon\to0^+$, the distribution of internal states remains regular where spatial gradients are order one, but may become singular where spatial gradients are small. }
    \label{fig:nontrivial_epsilon}
\end{figure}

To demonstrate the practical value of our results in the context of our specific model, we interpret our results in a biophysical context. In agreement with previous observations of phase separation in active systems \cite{auschra2021polarization,bialke2015active,buttinoni2013dynamical,cates2013brownian,fily2014freezing,fischer2019aggregation}, the non-uniform steady state profiles in Figure \ref{fig:bif_diagram_and_comparison} show that the final `phase-separated' state consists of a cluster of less motile (low $u$, low $\mathcal{D}$) cells coexisting with a dilute phase of more motile (high $u$, high $\mathcal{D}$) cells.  We plot an example `phase diagram' in Figure \ref{fig:phase_diagram} where we show the (normalised) variation in cell density $\Delta\rho/\rho^*$ of the stable phase separated state as a function of global cell density $\rho^*$ and the degree of motility repression $\mathcal{D}'_*$. Our weakly nonlinear analysis allows us to predict whether the phase transition is continuous or discontinuous, i.e.~whether $\Delta \rho$ varies continuously ($\mu>0$) or discontinuously ($\mu<0$) across the bifurcation point. 

Solutions on the non-trivial branch appear to remain regular as $\varepsilon\to0$, indicating that the instability leads to nontrivial chemical structuring within the population. This occurs despite the  terms in the perturbation expansion \eqref{eqn:wna_expand_n} being exponentially localized. We show the distribution of internal states in Figure \ref{fig:nontrivial_epsilon} as slices of $n(x,u)$ at fixed $x$.   The limit $\varepsilon\to0^+$ appears to be regular at $x=5.5$, but singular at $x=0.5$. Mathematically, this can be understood by noting that spatial gradients vanish on $x=0,L$ due to the no-flux conditions \eqref{eqn:gov_BCs_x}; therefore, $c(x)\approx C=\mathrm{const}.$ near the boundaries. As such, the governing equation \eqref{eqn:gov_eq_n} has an approximate solution near $x=0,L$ of the form \eqref{eqn:canonical_density_ss}, but with $u^*=g(C)/\lambda$. This means that near $x=0,L$, the non-uniform steady state $n(x,u)$ is sharply peaked (in $u$) for small $\varepsilon$.

\section{Discussion} \label{sec:discussion}

In summary, we have presented a WKBJ framework to extend classical weakly nonlinear analysis to structured models in settings where the base state is spatially uniform, but exponentially localized or singular with respect to the structured variables. As a specific example to demonstrate the effectiveness of our theory, we considered the model  of motile quorum sensing bacteria \eqref{eqn:gov_eqs}--\eqref{eqn:gov_BCs} from \cite{ridgway2023motility, ridgwaythesis}, reducing the local dynamics near the instability to the normal form of a pitchfork bifurcation \eqref{eqn:wna_C_ODE_simp}. From this normal form, we deduced the non-trivial steady states near the instability, as given in \eqref{eqn:wna_rho_nonuniform}. We additionally found that the quantity $\mu$ defined in \eqref{eqn:wna_mu_def} determines whether the bifurcation is subcritical ($\mu<0$) or supercritical ($\mu>0$). These results were compared with numerical steady state solutions and found to be in excellent agreement. This demonstrates the applicability of our framework.

Our analysis can be generalised in several ways. For example, we assumed that the reaction kinetics $f(u,c)$ and secretion rate $\alpha(u)$ were linear in $u$, while allowing for non-linearities in $c$. However, we could have considered fully nonlinear $f(u,c)$ with the mild constraint that there still exists a steady state solution (i.e.~a solution of \eqref{eqn:steady_algebraic}) and a few technical assumptions on the smoothness of $f$ and its growth at infinity.  In this case, there would be additional linear terms in the WKBJ amplitude equations \eqref{eqn:canonical_q_ODE}, \eqref{eqn:wna_adj_bif_condition}, and \eqref{eqn:wna_s_eqn}, and the phase of the exponentials in the perturbations in \eqref{eqn:wna_expand_n} would no longer be globally quadratic functions, but would remain locally quadratic. These complications are tractable within our framework as our analysis depends only on regular power series solutions to the amplitude equations and the use of Laplace's method to evaluate exponentially localized integrals. When $f$ is fully nonlinear, there are potentially multiple solutions of \eqref{eqn:steady_algebraic}, implying that the steady state is a sum of regularised Dirac deltas. In general, we would find that all but one of the Dirac masses have an exponentially small contribution to the steady state. As a future challenge, it may be interesting to explore the construction of quasi-steady solutions wherein the additional Dirac deltas have leading-order contributions to the linear and weakly nonlinear analyses.

As another potential generalisation, we expect that our analysis can be extended to consider Hopf bifurcations. As identified in \cite{ridgway2023motility}, a Hopf bifurcation occurs in \eqref{eqn:gov_eqs}--\eqref{eqn:gov_BCs} when motility is sufficiently promoted by quorum sensing. This condition is quantified by the condition $\mathcal{D}'(u^*)>\mathcal{D}_H'$ where $\mathcal{D}_H$ depends on the system parameters. Our framework would allow explicit determination of the limit cycle near the instability as well as quantification of when (and if) the bifurcation is subcritical.

Our framework can also be generalised to handle reaction terms in \eqref{eqn:gov_eq_n} as long as they are in balance for the spatially uniform steady state. For example, if we had included a logistic growth term on the RHS of \eqref{eqn:gov_eq_n} of the form $r(u)(1 - \rho/\rho_c)n$ for some constant carrying capacity $\rho_c$ and growth rate $r(u)$, then the linear instability criteria \eqref{eqn:eigenvalue_problem} for the non-trivial, spatially uniform steady state would become
\begin{equation}
    \sigma + D_ck^2 + \beta - \alpha_0\fc\rho_c\frac{\sigma + \left(\mathcal{D}_{0*} - u^*\mathcal{D}_*'\right)k^2 + r(u^*)}{\left(\sigma + \mathcal{D}_{0*}k^2 + r(u^*)\right)\left(\sigma + \mathcal{D}_{0*}k^2 + \lambda\right)}=\mathcal{O}(\varepsilon),
	\label{eqn:eigenvalue_problem_logistic}
\end{equation}
which can be obtained by a straightforward modification of our techniques. The weakly nonlinear analysis is a lengthy calculation but would be analytically tractable.

In terms of applications, this framework is  broadly applicable  to other structured models with sharply peaked solutions. One such example is  the stemness structured model of tumor heterogeneity in \cite{celora2023spatio}, but with modified boundary conditions. Since the inclusion of nonlocal terms is tractable within our framework, our work complements previous work on patterning in non-structured, nonlocal models \cite{carrillo2019aggregation,jewell2023patterning, topaz2006nonlocal, giunta2024weakly}.   Our techniques can also be applied to study systems where the internal dynamics arise through an advection term in state space, and where these dynamics have a steady state (akin to the point $u=u^*$ satisfying \eqref{eqn:steady_kinetics}). Internal dynamics will contribute an advection term whenever they can be described at the individual level by an ODE or system of ODEs (i.e.~a state space flow). This occurs frequently because many structured models can be `derived' through formal upscaling of individual-level models consisting of ODEs and/or SDEs. Our approach is therefore able to characterise the weakly nonlinear behavior of solutions in a wide range of structured models. 

\section*{Acknowledgments}

The authors would like to acknowledge the UCL Mathematics Department for use of the high performance cluster. P.P.~is supported by a UKRI Future Leaders Fellowship [MR/V022385/1]. W.J.M.R.~is supported by the Pacific Institute for the Mathematical Sciences. For the purposes of Open Access, the authors have applied a CC BY public copyright license to any Author Accepted Manuscript (AAM) version arising from this submission.

\bibliographystyle{siamplain}
\bibliography{bibliography_abbr}

\begin{thebibliography}{10}

\bibitem{al2022spikes}
{\sc F.~Al~Saadi, A.~Champneys, C.~Gai, and T.~Kolokolnikov}, {\em Spikes and
  localised patterns for a novel {S}chnakenberg model in the semi-strong
  interaction regime}, Eur. J. Appl. Math., 33 (2022), pp.~133--152.

\bibitem{arnold2012existence}
{\sc A.~Arnold, L.~Desvillettes, and C.~Pr{\'e}vost}, {\em Existence of
  nontrivial steady states for populations structured with respect to space and
  a continuous trait}, Commun. Pure Appl. Anal., 11 (2012), pp.~83--96.

\bibitem{auschra2021polarization}
{\sc S.~Auschra, V.~Holubec, N.~A. S{\"o}ker, F.~Cichos, and K.~Kroy}, {\em
  Polarization-density patterns of active particles in motility gradients},
  Phys. Rev. E, 103 (2021), p.~062601.

\bibitem{barles2009concentration}
{\sc G.~Barles, S.~Mirrahimi, and B.~Perthame}, {\em Concentration in
  {L}otka-{V}olterra {P}arabolic or {I}ntegral {E}quations: {A} {G}eneral
  {C}onvergence {R}esult}, Methods Appl. Anal., 16 (2009), pp.~321--340.

\bibitem{benichou2012front}
{\sc O.~B{\'e}nichou, V.~Calvez, N.~Meunier, and R.~Voituriez}, {\em Front
  acceleration by dynamic selection in fisher population waves}, Phys. Rev. E,
  86 (2012), p.~041908.

\bibitem{bialke2015active}
{\sc J.~Bialk{\'e}, T.~Speck, and H.~L{\"o}wen}, {\em Active colloidal
  suspensions: {C}lustering and phase behavior}, J. Non-Cryst. Solids, 407
  (2015), pp.~367--375.

\bibitem{buttinoni2013dynamical}
{\sc I.~Buttinoni, J.~Bialk{\'e}, F.~K{\"u}mmel, H.~L{\"o}wen, C.~Bechinger,
  and T.~Speck}, {\em Dynamical {C}lustering and {P}hase {S}eparation in
  {S}uspensions of {S}elf-{P}ropelled {C}olloidal {P}articles}, Phys. Rev.
  Lett., 110 (2013), p.~238301.

\bibitem{carrillo2019aggregation}
{\sc J.~A. Carrillo, K.~Craig, and Y.~Yao}, {\em Aggregation-{D}iffusion
  {E}quations: {D}ynamics, {A}symptotics, and {S}ingular {L}imits}, in Active
  Particles, Volume 2, N.~Bellomo, P.~Degond, and E.~Tadmor, eds., Springer,
  2019, pp.~65--108.

\bibitem{cates2013brownian}
{\sc M.~E. Cates and J.~Tailleur}, {\em When are active {B}rownian particles
  and run-and-tumble particles equivalent? {C}onsequences for motility-induced
  phase separation}, Europhys. {L}ett., 101 (2013), p.~20010.

\bibitem{cates2015motility}
{\sc M.~E. Cates and J.~Tailleur}, {\em Motility-{I}nduced {P}hase
  {S}eparation}, Annu. {R}ev. {C}ondens. {M}atter {P}hys., 6 (2015),
  pp.~219--244.

\bibitem{celora2023spatio}
{\sc G.~L. Celora, H.~M. Byrne, and P.~G. Kevrekidis}, {\em Spatio-temporal
  modelling of phenotypic heterogeneity in tumour tissues and its impact on
  radiotherapy treatment}, J. Theor. Biol., 556 (2023), p.~111248.

\bibitem{chapman1998exponential}
{\sc S.~Chapman, J.~King, and K.~Adams}, {\em Exponential asymptotics and
  {S}tokes lines in nonlinear ordinary differential equations}, Proc. R. Soc.
  Lon. Ser. A, 454 (1998), pp.~2733--2755.

\bibitem{chen2011stability}
{\sc W.~Chen and M.~J. Ward}, {\em The {S}tability and {D}ynamics of
  {L}ocalized {S}pot {P}atterns in the {T}wo-{D}imensional {G}ray--{S}cott
  {M}odel}, SIAM J. Appl. Dyn. Syst., 10 (2011), pp.~582--666.

\bibitem{crossley2025modelling}
{\sc R.~M. Crossley, P.~K. Maini, and R.~E. Baker}, {\em Modelling the {I}mpact
  of {P}henotypic {H}eterogeneity on {C}ell {M}igration: {A} {C}ontinuum
  {F}ramework {D}erived from {I}ndividual-{B}ased {P}rinciples}, Bull. Math.
  Biol., 87 (2025), p.~123.

\bibitem{cushing1985equilibria}
{\sc J.~Cushing}, {\em Equilibria in structured populations}, J. Math. Biol.,
  23 (1985), pp.~15--39.

\bibitem{cushing1998introduction}
{\sc J.~M. Cushing}, {\em An {I}ntroduction to {S}tructured {P}opulation
  {D}ynamics}, SIAM, Philadelphia, 1998.

\bibitem{dalwadi2021emergent}
{\sc M.~P. Dalwadi and P.~Pearce}, {\em Emergent robustness of bacterial quorum
  sensing in fluid flow}, Proc. Natl. Acad. Sci., 118 (2021), p.~e2022312118.

\bibitem{dalwadi2023universal}
{\sc M.~P. Dalwadi and P.~Pearce}, {\em Universal dynamics of biological
  pattern formation in spatio-temporal morphogen variations}, Proc. R. Soc. A,
  479 (2023), p.~20220829.

\bibitem{dewachter2019bacterial}
{\sc L.~Dewachter, M.~Fauvart, and J.~Michiels}, {\em Bacterial {H}eterogeneity
  and {A}ntibiotic {S}urvival: {U}nderstanding and {C}ombatting {P}ersistence
  and {H}eteroresistance}, Mol. {C}ell, 76 (2019), pp.~255--267.

\bibitem{diekmann2005dynamics}
{\sc O.~Diekmann, P.-E. Jabin, S.~Mischler, and B.~Perthame}, {\em The dynamics
  of adaptation: {A}n illuminating example and a {H}amilton--{J}acobi
  approach}, Theor. Popul. Biol., 67 (2005), pp.~257--271.

\bibitem{doelman1998stability}
{\sc A.~Doelman, R.~A. Gardner, and T.~J. Kaper}, {\em Stability analysis of
  singular patterns in the 1{D} {G}ray-{S}cott model: a matched asymptotics
  approach}, Physica D, 122 (1998), pp.~1--36.

\bibitem{doelman1997pattern}
{\sc A.~Doelman, T.~J. Kaper, and P.~A. Zegeling}, {\em Pattern formation in
  the one-dimensional {G}ray-{S}cott model}, Nonlinearity, 10 (1997), p.~523.

\bibitem{fily2014freezing}
{\sc Y.~Fily, S.~Henkes, and M.~C. Marchetti}, {\em Freezing and phase
  separation of self-propelled disks}, Soft Matter, 10 (2014), pp.~2132--2140.

\bibitem{fischer2019aggregation}
{\sc A.~Fischer, A.~Chatterjee, and T.~Speck}, {\em Aggregation and
  sedimentation of active {B}rownian particles at constant affinity}, J. Chem.
  Phys., 150 (2019), p.~064910.

\bibitem{genieys2006adaptive}
{\sc S.~G{\'e}nieys, V.~Volpert, and P.~Auger}, {\em Adaptive dynamics:
  modelling {D}arwin's divergence principle}, C.R. Biol., 329 (2006),
  pp.~876--879.

\bibitem{giunta2024weakly}
{\sc V.~Giunta, T.~Hillen, M.~A. Lewis, and J.~R. Potts}, {\em Weakly nonlinear
  analysis of a two-species non-local advection--diffusion system}, Nonlinear
  Anal. Real World Appl., 78 (2024), p.~104086.

\bibitem{gomez2021hopf}
{\sc D.~Gomez, L.~Mei, and J.~Wei}, {\em Hopf bifurcation from spike solutions
  for the weak coupling {G}ierer--{M}einhardt system}, Eur. J. Appl. Math., 32
  (2021), pp.~113--145.

\bibitem{heil2006oomph}
{\sc M.~Heil and A.~L. Hazel}, {\em oomph-lib--{A}n
  \emph{O}bject-\emph{O}riented \emph{M}ulti-\emph{Ph}ysics {F}inite-{E}lement
  \emph{Lib}rary}, in Fluid-Structure Interaction, H.-J. Bungartz and
  M.~Sch{\"a}fer, eds., vol.~53, Springer (Lecture Notes on Computational
  Science and Engineering), Berlin, Heidelberg, 2006, pp.~19--49.

\bibitem{hinch1991book}
{\sc E.~J. Hinch}, {\em Perturbation {M}ethods}, Cambridge University Press,
  Cambridge, 1991.

\bibitem{Hu2021}
{\sc R.~Hu, S.~Seager, and S.~Yasuda}, {\em Effects of internal dynamics on
  chemotactic aggregation of bacteria}, Phys. Biol., 18 (2021), p.~066001.

\bibitem{inaba2017age}
{\sc H.~Inaba}, {\em Age-{S}tructured {P}opulation {D}ynamics in {D}emography
  and {E}pidemiology}, Springer, Singapore, 2017.

\bibitem{iron2001stability}
{\sc D.~Iron, M.~J. Ward, and J.~Wei}, {\em The stability of spike solutions to
  the one-dimensional {G}ierer--{M}einhardt model}, Physica D, 150 (2001),
  pp.~25--62.

\bibitem{jabin2023collective}
{\sc P.-E. Jabin and B.~Perthame}, {\em Collective motion driven by nutrient
  consumption}, Asymptot. Anal., 133 (2023), pp.~483--497.

\bibitem{jewell2023patterning}
{\sc T.~J. Jewell, A.~L. Krause, P.~K. Maini, and E.~A. Gaffney}, {\em
  Patterning of nonlocal transport models in biology: {T}he impact of spatial
  dimension}, Math. Biosci., 366 (2023), p.~109093.

\bibitem{kolokolnikov2021competition}
{\sc T.~Kolokolnikov, F.~Paquin-Lefebvre, and M.~J. Ward}, {\em Competition
  instabilities of spike patterns for the 1{D} {G}ierer--{M}einhardt and
  {S}chnakenberg models are subcritical}, Nonlinearity, 34 (2021), p.~273.

\bibitem{kolokolnikov2009spot}
{\sc T.~Kolokolnikov, M.~J. Ward, and J.~Wei}, {\em Spot {S}elf-{R}eplication
  and {D}ynamics for the {S}chnakenburg {M}odel in a {T}wo-{D}imensional
  {D}omain}, J. Nonlinear Sci., 19 (2009), pp.~1--56.

\bibitem{kotbook}
{\sc M.~Kot}, {\em Elements of {M}athematical {E}cology}, Cambridge University
  Press, Cambridge, 2001.

\bibitem{krause2021modern}
{\sc A.~L. Krause, E.~A. Gaffney, P.~K. Maini, and V.~Klika}, {\em Modern
  perspectives on near-equilibrium analysis of {T}uring systems}, Philos.
  Trans. R. Soc. A, 379 (2021), p.~20200268.

\bibitem{krause2024pattern}
{\sc A.~L. Krause, V.~Klika, E.~Villar-Sep{\'u}lveda, A.~R. Champneys, and
  E.~A. Gaffney}, {\em Pattern {L}ocalization in the {S}wift--{H}ohenberg
  {E}quation via {S}lowly {V}arying {S}patial {H}eterogeneity}, SIAM J. Appl.
  Dyn. Syst., 24 (2025), pp.~2804--2847.

\bibitem{krause2020one}
{\sc A.~L. Krause, V.~Klika, T.~E. Woolley, and E.~A. Gaffney}, {\em From one
  pattern into another: analysis of turing patterns in heterogeneous domains
  via {WKBJ}}, J. R. Soc. Interface, 17 (2020), p.~20190621.

\bibitem{lorenzi2024phenotype}
{\sc T.~Lorenzi, N.~Loy, and C.~Villa}, {\em Phenotype-structuring of non-local
  kinetic models of cell migration driven by environmental sensing}, arXiv
  preprint arXiv:2412.16258,  (2024).

\bibitem{lorenzi2025derivation}
{\sc T.~Lorenzi, F.~R. Macfarlane, and K.~J. Painter}, {\em Derivation and
  travelling wave analysis of phenotype-structured haptotaxis models of cancer
  invasion}, Eur. J. Appl. Math., 36 (2025), pp.~231--263.

\bibitem{lorenzi2022trade}
{\sc T.~Lorenzi and K.~J. Painter}, {\em Trade-offs between chemotaxis and
  proliferation shape the phenotypic structuring of invading waves}, Int. J.
  Non-Linear Mech., 139 (2022), p.~103885.

\bibitem{lorenzi2025pattern}
{\sc T.~Lorenzi and K.~J. Painter}, {\em Pattern formation within
  phenotype-structured chemotactic populations}, arXiv preprint
  arXiv:2506.03389,  (2025).

\bibitem{lorenzi2025phenotype}
{\sc T.~Lorenzi, K.~J. Painter, and C.~Villa}, {\em Phenotype structuring in
  collective cell migration: a tutorial of mathematical models and methods}, J.
  Math. Biol., 90 (2025), p.~61.

\bibitem{lorenzi2022invasion}
{\sc T.~Lorenzi, B.~Perthame, and X.~Ruan}, {\em Invasion fronts and adaptive
  dynamics in a model for the growth of cell populations with heterogeneous
  mobility}, Eur. J. Appl. Math., 33 (2022), pp.~766--783.

\bibitem{lorenzi2020asymptotic}
{\sc T.~Lorenzi and C.~Pouchol}, {\em Asymptotic analysis of selection-mutation
  models in the presence of multiple fitness peaks}, Nonlinearity, 33 (2020),
  p.~5791.

\bibitem{lorz2011dirac}
{\sc A.~Lorz, S.~Mirrahimi, and B.~Perthame}, {\em Dirac mass dynamics in
  multidimensional nonlocal parabolic equations}, Commun. in Partial Differ.
  Equations, 36 (2011), pp.~1071--1098.

\bibitem{loy2024hamilton}
{\sc N.~Loy and B.~Perthame}, {\em A {H}amilton--{J}acobi approach to nonlocal
  kinetic equations}, Nonlinearity, 37 (2024), p.~105019.

\bibitem{auger2008structured}
{\sc P.~Magal and S.~Ruan}, {\em Structured {P}opulation {M}odels in {B}iology
  and {E}pidemiology}, Springer, Berlin, 2008.

\bibitem{marusyk2012intra}
{\sc A.~Marusyk, V.~Almendro, and K.~Polyak}, {\em Intra-tumour heterogeneity:
  a looking glass for cancer?}, Nat. Rev. {C}ancer, 12 (2012), pp.~323--334.

\bibitem{m1925applications}
{\sc A.~M'kendrick}, {\em Applications of {M}athematics to {M}edical
  {P}roblems}, Proc. Edinburgh Math. Soc., 44 (1925), pp.~98--130.

\bibitem{murraybook}
{\sc J.~D. Murray}, {\em Mathematical {B}iology. I, An {I}ntroduction},
  Interdisciplinary {A}pplied {M}athematics; 17, Springer, New York, 3rd
  ed.~ed., 2002.

\bibitem{olde1995stokes}
{\sc A.~Olde~Daalhuis, S.~Chapman, J.~King, J.~Ockendon, and R.~Tew}, {\em
  Stokes {P}henomenon and {M}atched {A}symptotic {E}xpansions}, SIAM J. Appl.
  Math., 55 (1995), pp.~1469--1483.

\bibitem{papenfort2016quorum}
{\sc K.~Papenfort and B.~L. Bassler}, {\em Quorum sensing signal--response
  systems in {G}ram-negative bacteria}, Nat. {R}ev. {M}icrobiol., 14 (2016),
  pp.~576--588.

\bibitem{perthame2007transport}
{\sc B.~Perthame}, {\em Transport {E}quations in {B}iology}, Springer,
  Birkh\"auser Basel, 2007.

\bibitem{perthame2007concentration}
{\sc B.~Perthame and S.~G{\'e}nieys}, {\em Concentration in the {N}onlocal
  {F}isher {E}quation: the {H}amilton-{J}acobi {L}imit}, Math. Modell. Nat.
  Phenom., 2 (2007), pp.~135--151.

\bibitem{phan2024direct}
{\sc T.~V. Phan, H.~H. Mattingly, L.~Vo, J.~S. Marvin, L.~L. Looger, and
  T.~Emonet}, {\em Direct measurement of dynamic attractant gradients reveals
  breakdown of the {P}atlak--{K}eller--{S}egel chemotaxis model}, Proc. Natl.
  Acad. Sci., 121 (2024), p.~e2309251121.

\bibitem{ridgwaythesis}
{\sc W.~J.~M. Ridgway}, {\em Motility-{I}nduced {P}atterning in {Q}uorum
  {S}ensing {B}acteria}, {DP}hil thesis, University of Oxford, 2025.

\bibitem{ridgway2023motility}
{\sc W.~J.~M. Ridgway, M.~P. Dalwadi, P.~Pearce, and S.~J. Chapman}, {\em
  Motility-{I}nduced {P}hase {S}eparation {M}ediated by {B}acterial {Q}uorum
  {S}ensing}, Phys. Rev. Lett., 131 (2023), p.~228302.

\bibitem{smith2023bacterial}
{\sc W.~P. Smith, B.~R. Wucher, C.~D. Nadell, and K.~R. Foster}, {\em Bacterial
  defences: mechanisms, evolution and antimicrobial resistance}, Nat. Rev.
  Microbiol., 21 (2023), pp.~519--534.

\bibitem{strogatz1994nonlinear}
{\sc S.~H. Strogatz}, {\em Nonlinear Dynamics and Chaos: With Applications to
  Physics, Biology, Chemistry, and Engineering}, Perseus Books, Reading, MA,
  1994.

\bibitem{topaz2006nonlocal}
{\sc C.~M. Topaz, A.~L. Bertozzi, and M.~A. Lewis}, {\em A {N}onlocal
  {C}ontinuum {M}odel for {B}iological {A}ggregation}, Bull. Math. Biol., 68
  (2006), pp.~1601--1623.

\bibitem{turanova2015model}
{\sc O.~Turanova}, {\em On a model of a population with variable motility},
  Math. Models Methods in Appl. Sci., 25 (2015), pp.~1961--2014.

\bibitem{turing}
{\sc A.~M. Turing}, {\em The chemical basis of morphogenesis}, Philos. Trans.
  R. Soc. Lon. Ser. B, 237 (1952), pp.~37--72.

\bibitem{veerman2015breathing}
{\sc F.~Veerman}, {\em Breathing pulses in singularly perturbed
  reaction-diffusion systems}, Nonlinearity, 28 (2015), p.~2211.

\bibitem{walker2010global}
{\sc C.~Walker}, {\em Global bifurcation of positive equilibria in nonlinear
  population models}, J. Differ. Equations, 248 (2010), pp.~1756--1776.

\bibitem{walker2025stability}
{\sc C.~Walker}, {\em Stability and {I}nstability of {E}quilibria in
  {A}ge-{S}tructured {D}iffusive {P}opulations}, J. Dyn. Differ. Equations, 37
  (2025), pp.~1315--1354.

\bibitem{wei2001spikes}
{\sc J.~Wei and M.~Winter}, {\em Spikes for the {T}wo-{D}imensional
  {G}ierer-{M}einhardt {S}ystem: {T}he {W}eak {C}oupling {C}ase}, J. Nonlinear
  Sci., 11 (2001), pp.~415--458.

\bibitem{wei2002spikes}
{\sc J.~Wei and M.~Winter}, {\em Spikes for {G}ierer-{M}einhardt system in
  {T}wo {D}imensions: {T}he {S}trong {C}oupling {C}ase}, J. Differ. Equations,
  178 (2002), pp.~478--518.

\bibitem{wong2020weakly}
{\sc T.~Wong and M.~J. Ward}, {\em Weakly {N}onlinear {A}nalysis of
  {P}eanut-{S}haped {D}eformations for {L}ocalized {S}pots of {S}ingularly
  {P}erturbed {R}eaction-{D}iffusion {S}ystems}, SIAM J. Appl. Dyn. Syst., 19
  (2020), pp.~2030--2058.

\end{thebibliography}

\appendix
\section{Regularity conditions for the WKBJ amplitudes} \label{app:stokes}

In this Appendix, we show that any solution of the WKBJ amplitude equations \eqref{eqn:canonical_q0_q1_ode}, \eqref{eqn:wna_wj_eqns}, \eqref{eqn:wna_sj_eqns} that satisfies the boundary conditions must be smooth at $u=u^*$. This effectively exchanges the boundary conditions with a local regularity condition, which we enforce in the main text through the use of a regular local power series. Since the only difference between \eqref{eqn:canonical_q0_q1_ode}, \eqref{eqn:wna_wj_eqns}, \eqref{eqn:wna_sj_eqns} is the inhomogeneous term on the RHS, the analysis is nearly identical for all three cases. We therefore only show the analysis for \eqref{eqn:canonical_q0_q1_ode}. Broadly, we demonstrate that $q_j(u)$ in \eqref{eqn:canonical_q0_q1_ode} must be smooth at $u=u^*$ by showing that any singularity in $q_j$ will switch on an exponentially large contribution in the far-field of $q$. If such a term were present, it would prevent $q$ from being polynomially bounded in the far-field as required by the boundary conditions. The switching of exponentially small terms in an asymptotic expansion is the mechanism behind Stokes phenomenon \cite{hinch1991book}, the apparent discontinuity of an asymptotic approximation of a continuous function. The main idea in Stokes phenomenon is that exponentially small terms can be switched on across certain curves (Stokes lines) in the complex plane. Even though these terms are buried deep in the series, they can become dominant across anti-Stokes lines, leading to an apparent discontinuity if the switching is not resolved by keeping track of the late terms in the series.

We begin by giving an outline of our methodology before getting into the technical details. First, we show that the `outer' expansion \eqref{eqn:canonical_q_expand} for $q(u)$ has a factorial-over-power divergence at large $j$ whenever there is a singular term in the series \eqref{eqn:canonical_q_expand}. To detect the Stokes phenomenon, we incorporate exponentially small terms into the expansion by optimally truncating the series and analysing the governing equation for the remainder. We then find that a term proportional to $\exp(\lambda(u-u^*)^2/2\varepsilon)$ is switched on across the Stokes line characterised by purely imaginary values of $u-u^*$. For purely real values of $u-u^*$, this term prevents  $q$ from being polynomially bounded as $u\to 0,\infty$, as required by the boundary conditions. Therefore, the only way to satisfy the boundary conditions is to demand that $q_j$ is smooth for all $j$, effectively exchanging the boundary conditions for a local regularity condition. Our analysis in this Appendix employs the methodology developed in \cite{chapman1998exponential,olde1995stokes} for detecting Stokes phenomenon in the solutions of linear and nonlinear ODEs.

Our first task is to show that if $q_0$ is non-smooth at $u^*$, then the late terms in the series \eqref{eqn:canonical_q_expand} exhibit a factorial-over-power divergence. To do this, we determine the large $j$ behaviour of $q_j$. Suppose that $q_j$ has a pole of order $p$ at $u^*$, then $q_{j+1}$ will have a pole of order $p+2$ at $u^*$. Intuitively, this is because $q_{j+1}$ is obtained from $q_j$ by differentiating $q_j$ twice, dividing by $u-u^*$ once, and integrating once. Motivated by this observation, we seek a large $j$ solution for $q_j$ in the form
\begin{equation}
	q_j(u)\sim Q\frac{\Gamma(j + \gamma(u) + 1)}{\lambda^j\big(v(u)\big)^{j + \kappa(u)}},
	\label{eqn:canonical_qj_large_j_ansatz}
\end{equation}
where $v$, $\gamma$, and $\kappa$ are functions independent of $j$ and $Q$ is a constant that will not be relevant in our analysis. The quantity $\Gamma$ is the Gamma function, which  can be expanded for large argument with a generalised form of Stirling's approximation (see e.g.~\cite{hinch1991book}, p.~34):
\begin{equation}
	\Gamma(z) = \sqrt{\frac{2\pi}{z}}\left(\frac{z}{e}\right)^z\left(1 + \frac{1}{12z} + \mathcal{O}(z^{-2})\right), \qquad \mathrm{Re}(z)>0.
	\label{eqn:stirling}
\end{equation}
To determine the unknown functions $v$, $\gamma$, and $\kappa$, we substitute \eqref{eqn:canonical_qj_large_j_ansatz} into \eqref{eqn:canonical_qj_ode} to obtain
\begin{align}
	-(sv'+(v')^2)j^3 &+ \left(sv + 2vv'\right)\gamma' j^2\log j - \bigg[\left(sv + 2vv'\right)\log v\kappa' +\left(sv' + 2(v')^2\right)\kappa \nonumber \\
	& -hv + sv'\left(2\gamma + \frac{25}{12}\right) + (v')^2\left(\gamma + \frac{11}{12}\right) - vv''\bigg]j^2 = o(j^2), \label{eqn:factorial_power_expand}
\end{align}
 where we define $s:=u-u^*$ and 
 \begin{equation}
 	h(u):=\frac{\mathcal{D}(u)k^2}{\lambda}>0,
 	\label{eqn:canonical_h_def}
 \end{equation}
 to simplify notation. Collecting terms at $\mathcal{O}(j^3)$ and requiring that $q_j$ has a pole at $s=0$, we have
\begin{equation}
	v(u) = -\frac{(u-u^*)^2}{2}.
	\label{eqn:stokes_v_sol}
\end{equation}
Then at next order, the $\mathcal{O}(j^2\log j)$ terms in \eqref{eqn:factorial_power_expand} yield $\gamma'=0$. We fix the constant value of $\gamma$ below. At next and final order, the $\mathcal{O}(j^2)$ terms in \eqref{eqn:factorial_power_expand} give the following ODE for $\kappa$
\begin{equation}
	s\kappa' \log v + 2\kappa  + h - 2\gamma  - 3 = 0.
    \label{eqn:kappa_ode}
\end{equation}
Eq.~\eqref{eqn:kappa_ode} has the general solution
\begin{equation}
	\kappa(u) = \gamma + \frac{3}{2} + \frac{\log q_h(u) + \tilde{C}}{\log\left(-\frac{(u-u^*)^2}{2}\right)},
	\label{eqn:kappa_sol}
\end{equation}
where $\tilde{C}$ is an integration constant, and $q_h$ is a homogeneous solution of \eqref{eqn:canonical_q0_q1_ode} given by
\begin{equation}
	q_h(u):= \exp\left(-\int_1^u\frac{h(\bar{u})}{\bar{u} - u^*}\,\mathrm{d}\bar{u}\right).
	\label{eqn:canonical_qh}
\end{equation}
The arbitrary lower bound of integration in \eqref{eqn:canonical_qh} is fixed for convenience. We set $\tilde{C}=0$ without loss of generality because any non-zero value could be absorbed into $Q$ in \eqref{eqn:canonical_qj_large_j_ansatz} due to the identity
\begin{equation}
	\left(-s^2/2\right)^{\tilde{C}/\log(-s^2/2)}=e^{\tilde{C}} = \mathrm{const}.
\end{equation}

The constant $\gamma$ is determined from demanding that the singularity structure of $q_j$ is consistent with the singularity structure of $q_0$. From \eqref{eqn:canonical_q0_q1_ode}, if $q_0\sim (u-u^*)^p$ as $u\to u^*$, then $q_j\sim (u-u^*)^{p + 2j}$.  To determine the strength of the singularity $p$ in $q_0$, we note that there is a particular solution $q_p$ of \eqref{eqn:canonical_q0_ode} that has the local behaviour 
\begin{equation}
	q_p(u) = \sqrt{\frac{\lambda^3}{2\pi}}\frac{\rho^*f_c^*C (u - u^*)}{\mathcal{D}(u^*)k^2 + \lambda} + \mathcal{O}\left((u-u^*)^2\right),\quad \text{as } u\to u^*.
	\label{eqn:canonical_qp_local}
\end{equation}
Since $q_p$ is smooth at $u^*$ the singularity structure of $q$ is determined by the poles of $q_h$. From \eqref{eqn:canonical_qh}, the local behaviour of $q_h$ is
\begin{equation}
	q_h(u)\sim a (u-u^*)^{-h_0}, \quad h_0:=h(u^*), \quad \text{as } u\to u^*,
	\label{eqn:canonical_qh_singular_structure}
\end{equation}
for some constant $a$. Thus $q_0$ has a pole of order $h_0$ and $q_j$ has a pole of order $h_0+2j$. Inserting  \eqref{eqn:kappa_sol} and \eqref{eqn:canonical_qh_singular_structure} into \eqref{eqn:canonical_qj_large_j_ansatz}, and demanding that the order of the pole for $q_j$ is equal to $h_0+2j$ gives
\begin{equation}
	\gamma = h_0 - \frac{3}{2}.
	\label{eqn:stokes_gamma_sol}
\end{equation}
We combine the large $j$ behaviour of $q_j$ in \eqref{eqn:canonical_qj_large_j_ansatz} with the expressions for $v$, $\kappa$, and $\gamma$ in Eqs.~\eqref{eqn:stokes_v_sol}, \eqref{eqn:kappa_sol}, and \eqref{eqn:stokes_gamma_sol} to find the asymptotic behaviour
\begin{equation}
	q_j(u) \sim Q\frac{(-1)^j\Gamma\left(j + h_0  - \frac{1}{2}\right)}{(u - u^*)^{2j + 2h_0}q_h(u)}\left(\frac{2}{\lambda}\right)^j \quad \text{as } j\to\infty.
	\label{eqn:canonical_qj_divergent}
\end{equation}
Thus the expansion \eqref{eqn:canonical_q_expand} for $q$ exhibits the expected factorial-over-power divergence. 

So far we have shown that if $q_0$ is non-smooth, then the expansion \eqref{eqn:canonical_q_expand} exhibits factorial-over-power divergence. We need to show that if the first non-smooth term  in the series  \eqref{eqn:canonical_q_expand} occurs at some $j=j'>0$, then \eqref{eqn:canonical_q_expand} still diverges as factorial-over-power.  For  $q_0,\ldots, q_{j'-1}$ to be smooth, there must be no contribution from the homogeneous solution $q_h$ until the term at $j=j'$. Then $q_{j'}$ (instead of $q_0)$ would have a pole of order $h_0$, implying a different value for the constant $\gamma$ which effectively  shifts the index $j$ in \eqref{eqn:canonical_qj_divergent}. So overall, the series \eqref{eqn:canonical_q_expand} still has a factorial-over-power divergence. Since the analysis that follows does not crucially depend on $\gamma$, we assume for simplicity that $q_0$ is non-smooth.

The next step is to incorporate exponentially small corrections into the expansion \eqref{eqn:canonical_q_expand} to detect Stokes switching. To do this, we introduce the remainder $R_N$ \emph{via}
\begin{equation}
	q(u) =\sum_{j=0}^{N-1} q_j(u)\varepsilon^j + R_N(u).
	\label{eqn:canonical_RN_def}
\end{equation}
Inserting \eqref{eqn:canonical_RN_def} into \eqref{eqn:canonical_q_ODE}, we find that $R_N$ exactly satisfies
\begin{equation}
	\varepsilon R_N'' - \lambda (u - u^*)R_N' - (\sigma + \mathcal{D}(u)k^2)R_N = -\varepsilon^N q_{N-1}''.
	\label{eqn:canonical_RN_ode}
\end{equation}
To determine the  number of terms $N$ in \eqref{eqn:canonical_RN_def}, we demand that the truncated series remain asymptotically valid. We therefore truncate the series when successive terms stop getting asymptotically smaller; mathematically this occurs for $j\geq N$, where $N$ satisfies the condition $\varepsilon^N |q_N| \sim \varepsilon^{N+1}|q_{N+1}|$. This gives $N\sim\lambda|u - u^*|^2/2\varepsilon$. Introducing complex polar coordinates as
\begin{equation}
	u - u^* = re^{i\theta}, \quad \theta \in [0,2\pi),
\end{equation} 
we fix $N$ as
\begin{equation}
	N=\frac{\lambda r^2}{2\varepsilon} + N_0(\varepsilon),
	\label{eqn:optimal_truncation}
\end{equation}
where $N_0=\mathcal{O}(1)$ as $\varepsilon\to0^+$. Since $N\to \infty$ as $\varepsilon\to0^+$, our truncation \eqref{eqn:optimal_truncation} is `optimal' in the sense that it incorporates all powers of $\varepsilon$. 

We detect Stokes phenomenon by determining the leading-order behaviour of $R_N(u)$ in different regions of the complex plane. We anticipate that a dominant contribution to the expansion \eqref{eqn:canonical_RN_def} from a homogeneous solution of \eqref{eqn:canonical_RN_ode} will be switched on across a Stokes line. Our next task therefore is to determine the leading-order (in $\varepsilon$) behaviour of the homogeneous solutions.  One such homogeneous solution is given by $R_N(u) \sim q_h(u)$ as $\varepsilon\to0^+$. Since the governing equation \eqref{eqn:canonical_q_ODE} for $q$ is linear, $q_h$ is a viable homogeneous solution at all orders of $\varepsilon$. As such, we absorb any contribution from $q_h$ to $R_N$ into the leading order $q_0$ term without loss of generality.  Since \eqref{eqn:canonical_RN_ode} is singularly perturbed, there is another solution which varies rapidly around $u^*$. We therefore seek a WKBJ solution in form 
\begin{equation}
	R_N(u) = \tilde{R}_N(u)\exp\left(\varepsilon^{-1}\psi(u)\right).
	\label{eqn:stokes_RN_wkbj}
\end{equation}
To determine the phase $\psi$, we insert \eqref{eqn:stokes_RN_wkbj} into \eqref{eqn:canonical_RN_ode} and collect the leading order powers of  $\varepsilon$ to obtain $\psi' = \lambda (u - u^*)$. Thus
\begin{equation}
	\psi(u) = \frac{\lambda}{2}(u-u^*)^2.
\end{equation}
The terms at next order in $\varepsilon$ yield
\begin{equation}
	\lambda (u - u^*)\tilde{R}_N' - \left(\sigma + \mathcal{D}k^2 - \lambda\right)\tilde{R}_N=0.
	\label{eqn:stoke_RN_amplitude_eqn}
\end{equation}
We solve \eqref{eqn:stoke_RN_amplitude_eqn} to find that the leading order behaviour of the homogeneous solutions of \eqref{eqn:canonical_RN_ode} are given by
\begin{align}
	 R_N&= \frac{1}{(u - u^*) q_h(u)}\left(1 + \mathcal{O}(\varepsilon)\right)\exp\left(\frac{\lambda (u - u^*)^2}{2\varepsilon}\right).
	\label{eqn:RN_hom_sol_large}
\end{align}

Since we expect a term proportional to \eqref{eqn:RN_hom_sol_large} to be switched on, we seek a particular solution of \eqref{eqn:stoke_RN_amplitude_eqn} in the form
\begin{equation}
	R_N(u) = \frac{S(u)}{(u - u^*)q_h(u)}\exp\left(\frac{\lambda (u- u^*)^2}{2\varepsilon}\right),
	\label{eqn:RN_wkbj_switch}
\end{equation}
where $S$ is the new dependent variable. We insert \eqref{eqn:RN_wkbj_switch} and \eqref{eqn:optimal_truncation} into \eqref{eqn:canonical_RN_ode}, and expand in $\varepsilon$, to find the following governing equation for  $S$
\begin{align}
	\dfo{S}{u} :=&  - \frac{i e^{-i\theta}}{r}\dfo{S}{\theta} \nonumber \\
	=& \frac{Q}{\varepsilon^{h_0}}\exp\left[-\frac{\lambda r^2}{2\varepsilon}(1 + e^{2i\theta}) - i\left(\frac{\lambda r^2}{2\varepsilon} + N_0\right)(2\theta - \pi) - 2ih_0\theta\right]\left[1  + \mathcal{O}(\varepsilon)\right],\label{eqn:stokes_switch_eqn}
\end{align}
where the constant prefactors that are independent of $N_0$ and $\varepsilon$ have been absorbed into $Q$. Since the RHS of \eqref{eqn:stokes_switch_eqn} is exponentially small in $\varepsilon$ everywhere except on $\theta = \pi/2$ and $\theta = 3\pi/2$, the function $S$ varies rapidly across the imaginary axis and is approximately constant everywhere else. This has the effect of switching on an exponentially small term proportional to $\exp(-\lambda r^2/2\varepsilon)$ across the imaginary axis, which defines the Stokes line. The contribution from this term becomes dominant across the anti-Stokes line $\mathrm{Re}((u-u^*)^2)=0$. To resolve the apparent discontinuity in the dominant terms of the expansion, we must explicitly resolve the switching in $S$ across the Stokes line. This allows us to connect the solutions on the real axis on either side of the singularity at $u=u^*$, which in turn allows us to impose the boundary conditions.

To resolve the Stokes switching, we consider a local inner region near $\theta=\pi/2$. With the benefit of hindsight, the appropriate scalings in the inner region are 
\begin{equation}
	\theta = \frac{\pi}{2} + \sqrt{\varepsilon} \tilde{\theta}, \qquad S = \frac{\tilde{S}}{\varepsilon^{h_0 - \frac{1}{2}}},
	\label{eqn:stokes_inner_var_def}
\end{equation}
where $\tilde{\theta}$ and $\tilde{S}$ are the inner variables. Inserting \eqref{eqn:stokes_inner_var_def} into \eqref{eqn:stokes_switch_eqn}, we find at leading order that $\tilde{S}$ satisfies
\begin{equation}
	\dfo{\tilde{S}}{\theta}\sim -\frac{Q r}{e^{i\pi(h_0 - \frac{1}{2})}}e^{-\lambda r^2\tilde{\theta}^2},
	\label{eqn:stokes_inner_S_eqn}
\end{equation}
Integrating \eqref{eqn:stokes_inner_S_eqn} and absorbing constants into $Q$, we obtain
\begin{equation}
	\tilde{S}(\tilde{\theta})\sim Q\left[\mathcal{S} - \mathrm{erf}\left(\sqrt{\lambda} r \tilde{\theta}\right)\right],
	\label{eqn:canonical_inner_switch}
\end{equation}
where  $\mathrm{erf}(\cdot)$ is the error function and $\mathcal{S}$ is a constant. 

The final step is to show that $q$ cannot satisfy the boundary conditions for any value of the constant $\mathcal{S}$ in \eqref{eqn:canonical_inner_switch}. To this end, we insert the solution for $\tilde{S}$ in \eqref{eqn:canonical_inner_switch}, along with the scalings \eqref{eqn:stokes_inner_var_def}, into $R_N$ in \eqref{eqn:RN_wkbj_switch} to find the leading order behaviour
\begin{equation}
	R_N\sim \frac{Q\varepsilon^{\frac{1}{2} - h_0}}{(u - u^*)q_h(u)}\exp\left[{\frac{\lambda (u - u^*)^2}{2\varepsilon}}\right]\left[\mathcal{S} - \mathrm{erf}\left(\sqrt{\frac{\lambda}{\varepsilon}}r\left(\frac{\pi}{2} - \theta\right)\right)\right], \label{eqn:canonical_RN_full_sol} 
\end{equation}
where $u - u^*:=re^{i\theta}$. Since $R_N$ must be polynomially bounded  for $q$ to be polynomially bounded (for real $u$), the terms in square brackets in \eqref{eqn:canonical_RN_full_sol} must  be exponentially small to balance the exponentially large factor. This is impossible to achieve with only one free constant $\mathcal{S}$. The only way to satisfy the boundary conditions is to take $Q=0$, which amounts to demanding that none of the terms in the expansion \eqref{eqn:canonical_q_expand} are singular at $u=u^*$. Any singular term will set off a factorial-over-power divergence, which in turn leads to Stokes phenomenon wherein  a contribution $R_N$ is switched on, again making it impossible to satisfy the boundary conditions.

\end{document}